\documentclass[10pt,reqno]{amsart}
\usepackage{fullpage}
\usepackage{times}
\usepackage{graphicx}

\usepackage{amsmath,amsthm,amsfonts,amscd,amssymb}
\usepackage{pst-node}
\usepackage{pst-plot}

\usepackage{a4}
\usepackage{amssymb}
\usepackage{array}
\usepackage{booktabs}
\usepackage{multirow}
\usepackage{hhline}

\newtheorem{thm}{Theorem}

\newtheorem{cor}[thm]{Corollary}
\newtheorem{lem}[thm]{Lemma}
\newtheorem{prop}[thm]{Proposition}
\newtheorem{defn}[thm]{Definition}
\newtheorem{ex}[thm]{Example}

\newtheorem{rem}[thm]{Remark}

\numberwithin{thm}{section}
\numberwithin{equation}{section}
\newcommand{\Real}{\mathbb R}

\newcommand{\eps}{\varepsilon}

\newcommand{\la}{\langle}
\newcommand{\ra}{\rangle}

\newcommand{\Comp}{\mathbb{C}}

\newcommand{\Hi}{\mathcal{H}}

\newcommand{\n}{\mathbb{N}}

\newcommand{\tor}{\mathbb{T}}

\newcommand{\z}{\mathbb{Z}}

\newcommand{\om}{\omega}

\newcommand{\prt}{\widehat{\otimes}}

\newcommand{\fM}{{\mathbf{M}}}

\begin{document}
\title{Twisted Fourier(-Stieltjes) spaces and amenability}

\author{Hun Hee Lee}
\author{Xiao Xiong}

\address{Hun Hee Lee : Department of Mathematical Sciences and Research Institute of Mathematics, Seoul National University, Gwanak-ro 1, Gwanak-gu, Seoul, 08826, Republic of Korea}
\email{hunheelee@snu.ac.kr}

\address{Xiao Xiong : Institute for Advanced Study in Mathematics, Harbin Institute of Technology, 150001 Harbin, China}
\email{xxiong@hit.edu.cn}


\thanks{HHLee was supported by the Basic Science Research Program through the National Research Foundation of Korea (NRF) Grant RS-2022-NR069971.
XXiong was supported by the National Natural Science Foundation of China (Grant No. 12371138 and No. W2441002).} 

\begin{abstract}
The Fourier(-Stieltjes) algebras on locally compact groups are important commutative Banach algebras in abstract harmonic analysis. In this paper we introduce a generalization of the above two algebras via twisting with respect to 2-cocycles on the group. We also define and investigate basic properties of the associated multiplier spaces with respect to a pair of 2-cocycles.
We finally prove a twisted version of the results of Nebbia, Fendler, Bo\.{z}ejko, Losert and Ruan characterizing amenability of the underlying locally compact group through comparison of the twisted Fourier-Stieltjes space with the associated multiplier spaces.
\end{abstract}

\maketitle


\section{Introduction}

One popular theme in abstract harmonic analysis is to construct Banach or operator algebras out of locally compact groups reflecting the properties of the underlying group. The first such examples would be the group algebra $L^1(G)$ and the measure algebra $M(G)$ for a locally compact group $G$ with respect to the convolution product. 
Some group properties are well represented by these algebras. For example, the group algebra $L^1(G)$ is amenable as a Banach algebra in the sense of B.E. Johnson if and only if $G$ is amenable as a locally compact group \cite{Johnson72} if and only if $M(G)$ is Connes amenable \cite{Runde03}.
Left translation on $G$ provides an important example of unitary representations of $G$, namely the left regular representation $\lambda: G \to \mathcal{U}(L^2(G))$, whose canonical lifting to the $L^1(G)$ leads us to the reduced group $C^*$-algebra $C^*_r(G)$ and the group von Neumann algebra $VN(G)$ of $G$ by completing the image of the lifting in the norm topology and in the strong operator topology in $B(L^2(G))$, respectively.
For discrete group $G$ it is well-known that $G$ is amenable if and only if $C^*_r(G)$ is nuclear as a $C^*$-algebra if and only if $VN(G)$ is semidiscrete as a von Neumann algebra \cite[Section 2.6]{BO08}.

While $L^1(G)$ and $M(G)$ are relatively easy to understand as Banach spaces, they are non-commutative Banach algebras (without $C^*$-algebra structures) for non-abelian groups, which makes the analysis for their algebraic structures quite limited.
This limitation leads us to look for other alternatives, and we do have dual objects which are commutative Banach algebras, namely the Fourier algebra $A(G)$ and the Fourier-Stieltjes algebra $B(G)$. 
It is well known that $A(G)^* \cong VN(G)$ and $B(G) \cong (C^*(G))^*$ in canonical ways, where $C^*(G)$ is the full $C^*$-algebra of $G$, which is the enveloping $C^*$-algebra of the Banach $*$-algebra $L^1(G)$ with respect to the canonical involution.
See e.g. \cite{EYM1964,KanLau} for more details for these algebras and the dual relations; a brief description is also given in Section \ref{subsec-fs-oa}.
The duality $A(G)^* \cong VN(G)$ allows us to consider a canonical operator space structure on $A(G)$ and the following result of Ruan \cite{Ruan95}: the algebra $A(G)$ is amenable as a completely contractive Banach algebra if and only if $G$ is amenable as a locally compact group.
See Section \ref{subsec-OS} for the details on operator spaces.
In this dual side, we actually have a richer structure, namely the algebra $MA(G)$ of bounded multipliers of $A(G)$ and the algebra $M_{cb}A(G)$ of completely bounded multipliers of $A(G)$, both containing $B(G)$. Note that for the corresponding objects $ML_1(G)$ and $M_{cb}L_1(G)$ in the group algebra side we have $ML^1(G)=M_{cb}L^1(G)=M(G)$ \cite{Wen52}. This additional structure gives us another way of seeing the underlying group structure.
The result of Nebbia \cite{Nebbia} says that, for a discrete group $G$, the condition $B(G) = MA(G)$ characterizes amenability of $G$.
The corresponding result for $M_{cb}A(G)$ was proved by Fendler in an unpublished note.
A closely related result was given by Bo\.{z}ejko \cite{Boz85}, which we will revisit in Section \ref{subsec-discrete}.
For a non-discrete locally compact group $G$ the equivalence between the condition $B(G) = MA(G)$ and the amenability of $G$ has been established by Losert \cite{Losert} and the corresponding result for $M_{cb}A(G)$ was proved by Ruan \cite{Ruan-unpub}, which is also unpublished.
We will carefully examine the approach of the above two results in this paper.

The above line of research has been extended to the cases of locally compact quantum groups (see \cite{CLR15} and \cite{Ruan96}). While the class of locally compact quantum groups is a vastly bigger one than the class of locally compact groups, it misses some of the fundamentally important examples of ``{\it quantum spaces}'' such as non-commutative tori (or quantum tori) and quantum Euclidean spaces (or Moyal planes). Recall that ``quantum spaces'' including locally compact quantum groups are understood through the ``algebras of functions on them'', which are either $C^*$-algebras or von Neumann algebras.
For example, the ``algebra of continuous functions'' on the non-commutative tori is the well-known $C^*$-algebra, called the irrational rotation algebra $C(\tor^2_\theta)$. Note that $C(\tor^2_\theta)$ does not allow any (locally) compact quantum group structure on it due to a result of P. Soltan \cite{Sol10}.
Nevertheless, there are canonical ways of obtaining $C(\tor^2_\theta)$ from locally compact groups, namely via 2-cocyle twisting.
More precisely, we begin with the dual group $\z^2 = \widehat{\tor^2}$ and a 2-cocycle
$$\sigma: \z^2 \times \z^2 \to \tor,\; ((x,y), (x',y')) \mapsto e^{i\theta xy'}$$
for a fixed irrational number $\theta$.
Then, the space $\ell^1(\z^2)$ carries a Banach $*$-algebra structure with respect to the twisted convolution $*_\sigma$ and the twisted involution $\star$, which will be clarified later. Now we get $C(\tor^2_\theta)$ as its enveloping $C^*$-algebra. This procedure of 2-cocyle twisting can be done for general locally compact groups and that is one of the main focuses of this paper. Note that the study of locally compact groups with 2-cocycle twisting goes back to the results of Mackey \cite{Mackey}, Edwards/Lewis \cite{EL} and Kleppner/Lipsman \cite{KL72}.

As is demonstrated above, 2-cocyles on locally compact groups give us more of Banach/operator algebras, namely twisted group/measure algebras and twisted group $C^*$-algebras and group von Neumann algebras,
which are different from the (quantum) group cases, and they have been studied extensively (see \cite{Packer} and the continued works for example).
However, analogues of the Fourier(-Stieltjes) algebra and the associated multiplier algebra received relatively little attention until now, except the case of discrete groups by B\'edos/Conti \cite{BC2009, BC2016-2}. It is our main goal to initiate a systematic study on the twisted Fourier(-Stieltjes) spaces and their associated multiplier spaces for any locally compact group. Note that we use the expression ``Fourier(-Stieltjes) spaces'' instead of ``Fourier(-Stieltjes) algebras'' since they do not carry Banach algebra structures with respect to the canonical operations in general. 
Instead, they carry bimodule structures with respect to the usual (untwisted) Fourier(-Stieltjes) algebras of the same underlying group. This fundamental difference leads us to consider multiplier spaces between the usual Fourier algebra and the twisted Fourier space (see Section \ref{sec-twisted-multiplier} for the details), 
which turn out to be the right object to establish a twisted version of the results of Nebbia, Fendler, Bo\.{z}ejko, Losert and Ruan \cite{Nebbia, Boz85, Losert, Ruan-unpub} characterizing amenability through the associated multiplier spaces.

\medskip

\textbf{Organization.} \color{black}
This paper is organized as follows. In Section \ref{sec-prelim} we collect basic materials needed in the sequel, namely various function spaces/algebras and operator algebras on locally compact groups.
We also cover necessary preliminaries on operator spaces.
In Section \ref{sec-twisted-alg} we recall known constructions of several twisted algebras on locally compact groups with more details scattered in the literature.
In Section \ref{sec-twisted-space} we define our main objects, namely twisted Fourier(-Stieltjes) spaces, and develop basic theory on them in parallel with the untwisted case. In Section \ref{sec-twisted-multiplier} we continue to define twisted multiplier spaces with a characterization in the style of Gilbert \cite{Gil}  and Jolissaint \cite{Joli92}.
In the final section, we address our main result, which is a characterization of amenability through twisted multiplier spaces after Nebbia, Fendler, Bo\.{z}ejko, Losert and Ruan.
The proof is divided into two sub-cases, namely the discrete group case (Section \ref{subsec-discrete}) and the non-discrete group case (Section \ref{subsec-non-discrete}).
Note that the original proof of Losert \cite{Losert} for the latter case contains an error, which remained unnoticed by the abstract harmonic analysis community even after a polished and more detailed version was presented in the recent monograph \cite[Proposition 5.3.3]{KanLau}.
Moreover, its operator space version by Ruan \cite{Ruan-unpub} was never published. See Remark \ref{rem-comments-proof-Losert} and \ref{rem-comment-cb} for the details.

\section{Preliminaries}\label{sec-prelim}

\subsection{Related function spaces and operator algebras without twisting}\label{subsec-fs-oa}

In this subsection we will collect basic materials regarding relevant function spaces and operator algebras associated to a locally compact group $G$.
We consider {\it the group algebra}, which is the space $L^1(G)$ with respect to a fixed left Haar measure $ds$ on $G$ equipped with the convolution $*$ and the involution $^*$ given by $$f*g(s) := \int_G f(t) g(t^{-1}s)dt,\;\; f^*  (s) = \Delta_G(s^{-1})\,\overline{f(s^{-1})},\;\; f,g\in L^1(G).$$
Here, $\Delta_G$ is the {\it modular function} of $G$. The enveloping $C^*$-algebra $C^*(G)$ of the Banach $*$-algebra $(L^1(G), *)$ is called {\it the full group $C^*$-algebra} of $G$. The enveloping von Neumann algebra $W^*(G)$ of $(L^1(G), *)$ (or of $C^*(G)$) can be identified with $(C^*(G))^{**}$ in a canonical way. There are reduced versions of the above algebras.
We begin with {\it the left regular representation} 
$$\lambda: G \to \mathcal{U}(L^2(G)),\;\lambda(s)f(t) = f(s^{-1}t),\;f\in L^2(G),\; s,t \in G,$$
which can be lifted to $L^1(G)$ by 
$$\tilde{\lambda}:L^1(G) \to B(L^2(G)),\; f \mapsto \tilde{\lambda}(f) := \int_G f(s)\lambda(s)ds.$$
Then we get {\it the reduced $C^*$-algebra} $C^*_r(G) := \overline{\tilde{\lambda}(L^1(G))}\subseteq B(L^2(G))$ and {\it the group von Neumann algebra} $VN(G) := \tilde{\lambda}(L^1(G))^{''}\subseteq B(L^2(G))$.
The dual of the full group $C^*$-algebra $(C^*(G))^*$ can be identified with {\it the Fourier-Stieltjes algebra} 
$$B(G) := \{\la \pi(\cdot)\xi, \eta \ra|\, \pi: G \to \mathcal{U}(\Hi_\pi)\, \text{unitary representation},\; \xi,\eta\in \Hi_\pi\},$$
which can be obtained by the linear span of $\mathcal{P}_1(G)$, the collection of all bounded {\it positive definite} functions on $G$ with value 1 at the identity $e$ of $G$.
The group von Neumann algebra $VN(G)$ is known to have a unique predual $A(G)$ called {\it the Fourier algebra}, which can be canonically identified as a subalgebra of $C_0(G)$.
Thanks to the fact that $VN(G) \subseteq B(L^2(G))$ is a standard form (\cite{Tak2}) we can deduce that
$$A(G) = \{\la \lambda(\cdot)f, g \ra|\, f,g\in L^2(G)\},$$
i.e. the space of coefficient functions of the left regular representation.
Both of $A(G)$ and $B(G)$ are Banach algebras with respect to the pointwise multiplication and $A(G)$ actually is an ideal of $B(G)$.

We will need one more function algebra $LUC(G)$ consisting of bounded continuous functions $f: G\to \Comp$ for which $s \in G \mapsto \lambda(s)f$ is continuous with respect to $\|\cdot\|_{\infty}$, which is clearly a commutative unital $C^*$-algebra.

{\it Amenability} is one of the most important properties of a locally compact group. We say that $G$ is amenable if the commutative von Neumann algebra $L^ \infty(G)$ has a left invariant mean $m$, i.e. a state on $L^ \infty(G)$ such that $m(\lambda(s)f)=m(f),~\forall s\in G,~f\in  L^ \infty(G).$ It is well-known that $G$ is amenable if and only if the canonical quotient map $Q: C^*(G) \to C^*_r(G)$ is injective \cite{Hul64}.

\subsection{Operator spaces}\label{subsec-OS}
A {\it concrete operator space} is a norm-closed subspace $E$ of $B(H)$ for some Hilbert space $H$. Immediate examples of operator spaces are, of course, $C^*$-algebras and von Neumann algebras. The following materials can be found in \cite{ER00}.

For an operator space $E\subseteq B(H)$ we can associate a sequence of Banach spaces $M_n(E)\subseteq M_n(B(H)) \cong B(\ell^2_n \otimes_2 H)$, $n\in \n$, whose canonical norms are denoted by $\|\cdot\|_n$.
Here, $\otimes_2$ refers to the tensor product of Hilbert spaces.
This sequence of norms $(\|\cdot\|_n)_{n\ge 1}$ (which we call an {\it operator space structure}) satisfies the {\it Ruan's axioms}:
	$$(R_1)\;\; \|\alpha x \beta\|_n \le \|\alpha\| \cdot \|x\|_n \cdot \|\beta\|,\;\; \text{and}\;\; (R_2)\;\; \|x\oplus y\|_{n+m} = \max\{ \|x\|_n, \|y\|_m \}$$
for any $\alpha, \beta \in M_n$, $x \in M_n(E)$ and $y\in M_m(E)$, $n,m\in \n$.
It turns out that $(R_1)$ and $(R_2)$ actually characterize operator spaces in the sense that any Banach space $E$ equipped with a sequence of norms $\|\cdot\|_n$ on $M_n(E)$ satisfying $(R_1)$ and $(R_2)$ can be isometrically embedded in $B(H)$ for some Hilbert space $H$.

For operator spaces $E$ and $F$ we say that a linear map $T:E\to F$ is {\it completely bounded} (shortly, a {\it cb-map}) if the {\it cb-norm} 
$$\|T\|_{cb} := \sup_{n\ge 1}\|I_n \otimes T: M_n(E) \to M_n(F)\|$$
is finite.
The space of all cb-maps from $E$ to $F$ will be denoted by $CB(E,F)$ equipped with the cb-norm. The map $T$ is called a {\it complete isometry} if $I_n \otimes T: M_n(E) \to M_n(F)$ is an isometry for all $n\ge 1$. It is called a {\it complete quotient map} if $I_n \otimes T: M_n(E) \to M_n(F)$ is a quotient map for all $n\ge 1$.

Ruan's axioms allow us to consider a natural duality for operator spaces since we can equip $M_n(E^*)$ with the norm  $\|\cdot \|'_n$coming from $CB(E,M_n)$, which can be shown to satisfy Ruan's axioms again. We call $(E^*, (\|\cdot \|'_n)_{n\in \n})$ the {\it operator space dual} of $E$. If the operator space $E$ is a dual Banach space, so that $E = (E_*)^*$ for some predual $E_*$ of $E$, then the usual isometric embedding $E_* \hookrightarrow  (E_*)^{**} = E^*$ provides a canonical operator space structure on $E_*$. For a linear map $T:E\to F$ between operator spaces, it is well-known that $T$ is a complete quotient map if and only if $T^*:E^*\to F^*$ is a complete isometry.

Among various tensor products available for operator spaces the {\it projective operator space tensor product} $E\prt F$ of two operator spaces $E$ and $F$ are characterized by the following property: the canonical map gives a complete isometry $(E \prt F)^* \cong CB(E,F^*)$. For a dual operator space $E = (E_*)^*$ we have $M_n(E)_* \cong (M_n)_* \prt E_*$ completely isometrically. Moreover, for two complete quotient maps $T_1:E_1 \to F_1$ and $T_2: E_2\to F_2$ between operator spaces their tensor product $T_1\otimes T_2: E_1 \prt E_2 \to F_2 \prt F_2$ is also a complete quotient map.
 
For a Hilbert space $H$ we can associate two important operator spaces $H_R := B(H,\Comp)\subseteq B(H)$ and $H_C := B(\Comp, H)\subseteq B(H)$ called the {\it row Hilbert space} and the {\it column Hilbert space}, respectively. When $H = \ell^2_n$, we simply write $R_n$ and $C_n$, respectively. It is well known that $H_R \prt H_C \cong B(H)_*$ completely isometrically via the canonical map. In particular, we have $R_n \prt C_n \cong (M_n)_*$.
For another Hilbert space $K$ we have complete isometries 
$$H_C \prt K_C \cong (H\otimes_2 K)_C,\;\; H_R \prt K_R \cong (H\otimes_2 K)_R.$$

\section{Twisted algebras on locally compact groups}\label{sec-twisted-alg}

In this section we present known constructions of several algebras with 2-cocycle twisting on a locally compact group $G$. We will always assume the second countability of $G$, which is necessary in the properties \eqref{eq-Borel-embed} and \eqref{eq-local-conti-embed} below. We begin with recalling the definition of 2-cocycles.

\begin{defn}
A Borel measurable function $\sigma : G\times G \to \tor$ is called a {\bf 2-cocycle} on $G$ if
$$\begin{cases}\sigma(s,t)\sigma(st,r) = \sigma(s,tr)\sigma(t,r),\\ \sigma(s,e) = \sigma(e,t) = 1 \end{cases}$$
for any $s,t,r\in G.$
The collection of all 2-cocycles on $G$ will be denoted by $Z^2(G,\tor)$, which is a group under pointwise multiplication.
\end{defn}
	
\begin{ex}\label{ex-2-cocyles}
Let $A \in M_n(\Real)$. Then,
\begin{equation}\label{eq-2-cocyles-ex}
\sigma_A: \Real^n \times \Real^n \to \tor,\;\;  (s,t) \mapsto  e^{i\la s, At \ra}
\end{equation}
is a 2-cocycle on the group $\Real^n$.
The same function gives rise to a 2-cocyle on the group $\z^n$ by restriction, which we still denote by $\sigma_A$ by abuse of notation.
For more examples of 2-cocycles, see \cite{Packer} and \cite{BCLPY}, where the latter contains 2-cocycles with discontinuity.
\end{ex}

\begin{defn}
Two 2-cocycles $\sigma_1, \sigma_2 \in Z^2(G,\tor)$ are called {\bf similar} (or {\bf equivalent}) if there is a Borel measurable function $\xi: G \to \tor$ such that
$$\sigma_1(s,t) = \frac{\xi(s)\xi(t)}{\xi(st)}\sigma_2(s,t),\;\; s,t\in G$$
satisfying $\xi(e) = 1$.
\end{defn}
		
Any 2-cocycle $\tau \in Z^2(G,\tor)$ is similar to a 2-cocycle $\sigma \in Z^2(G,\tor)$ satisfying the following normalization condition (\cite[p.215]{Kleppner74}):
\begin{equation}\label{eq-add-normalization}
\sigma(s,s^{-1}) = 1,\;\; \forall s\in G.
\end{equation}
Indeed, we can choose the map
$$\sigma(s,t) = \tau(st, (st)^{-1})^{\frac{1}{2}}\tau(s, s^{-1})^{-\frac{1}{2}}\tau(t, t^{-1})^{-\frac{1}{2}}\tau(s,t),\;\; s,t\in G.$$
It is straightforward to see that the above condition \eqref{eq-add-normalization} implies that
\begin{equation}\label{eq-add-normalization-consequence}
\sigma(s,t) = \overline{\sigma(t^{-1}, s^{-1})},\;\; s,t\in G.
\end{equation}
Combined with the above remark we may assume that the given 2-cocycle $\sigma \in Z^2(G,\tor)$ satisfies the above condition \eqref{eq-add-normalization} when we are only interested in the $*$-algebra structure of $L^1(G,\sigma)$.
Along with this discussion we say that $\sigma \in Z^2(G,\tor)$ is {\bf trivial} if it is similar to the constant 1 function, which obviously is a 2-cocycle.
				
Furthermore, any 2-cocycle $\sigma \in Z^2(G,\tor)$ is known (\cite[Theorem 3.1]{Par}) to be similar to a 2-cocycle $\tau \in Z^2(G,\tor)$ which is locally continuous at the identity, i.e. there is an open neighborhood $U$ of the identity of $G$ such that $\tau$ is continuous on $U\times U$.
For this reason, we will assume that $\sigma \in Z^2(G,\tor)$ is locally continuous at the identity from now on.
We also note that, even though we will not need \eqref{eq-add-normalization} or \eqref{eq-add-normalization-consequence} in this paper,
any 2-cocycle $\sigma \in Z^2(G,\tor)$ is similar to a 2-cocycle $\tau \in Z^2(G,\tor)$ which is both normalized and locally continuous at the identity, since the construction in \cite[Theorem 3.1]{Par} keeps the normalization condition \eqref{eq-add-normalization}.

The nature of ``twisting" coming from a 2-cocyle $\sigma \in Z^2(G,\tor)$ can be observed through the following {\bf central extension} of $G$.
Let $G_\sigma = G\times \tor$ (as a set) endowed with the group law
$$(s,z) \cdot (t,w) = (st, zw\sigma(s,t)),\;\; s,t\in G,\; z,w\in \tor.$$
The group $G_\sigma$ can be equipped with the Weil topology, which makes it a second countable locally compact group.
Recall that the Weil topology comes from the neighborhood system of the identity consisting of the sets $EE^{-1}$, where $E\subseteq G_\sigma$ is a Borel set of finite positive measure with respect to $dsdz$.
Here, $ds$ is a fixed left Haar measure on $G$ and $dz$ is the normalized Haar measure on $\tor$. Note that $dsdz$ is a left Haar measure on $G_\sigma$, which we will fix from now on. 

Thanks to second countability of $G$ we know that the Borel structure on $G_\sigma$ is the same as the product Borel structure on $G\times \tor$ (see \cite[p.13]{Par}),
which means that the formal identity $G\times \tor \to G_\sigma$ is a Borel isomorphism.
In particular,
\begin{equation}\label{eq-Borel-embed}
\text{the map $s\in G \mapsto (s,1)\in G_\sigma$ is Borel measurable.}    
\end{equation}
Note that the Weil topology on $G_\sigma$ and the product topology on $G\times \tor$ do not coincide unless $\sigma$ is continuous.
Since we are assuming $\sigma$ is locally continuous at the identity,
the Weil topology on $G_\sigma$ and the product topology on $G\times \tor$ coincide on a neighborhood of the identity,
so that
\begin{equation}\label{eq-local-conti-embed}
\text{the map $s\in G \mapsto (s,1)\in G_\sigma$ is also locally continuous at the identity.}  
\end{equation}

\subsection{Twisted convolution algebras}

The space $L^1(G)$ can be equipped with a Banach $*$-algebra structure different from the usual group algebra structure as follows.
\begin{defn}
Let $\sigma \in Z^2(G,\tor)$. For $f,g \in L^1(G)$ we define the {\bf twisted convolution} $f*_\sigma g$ and {\bf twisted involution} $f^\star$ by
$$\begin{cases}\displaystyle f*_\sigma g(s) := \int_G f(t) \sigma(t,t^{-1}s)g(t^{-1}s)dt,\\ \displaystyle f^\star (s) := \Delta_G(s^{-1})\overline{\sigma(s,s^{-1})}\,\overline{f(s^{-1})}\end{cases}$$
	for $s\in G$. Here, $\Delta_G$ is the modular function of $G$. The Banach $*$-algebra $(L^1(G), *_\sigma, \star)$ will be simply denoted by $L^1(G,\sigma)$.
	\end{defn}

\begin{rem}
\begin{enumerate}
\item It is immediate to see that the twisted convolution is associative on $L^1(G)$.
\item Since we have $|f*_\sigma g(s)|\le | \, |f| * |g|(s)\, |$ for $f, g\in C_c(G)$, $s\in G$ we can easily see that $f*_\sigma g$ is well-defined for $f,g\in L^1(G)$ as an $L^1$-function on $G$. 
Moreover, Young's inequality for the usual convolution tells us that its twisted version also holds. In other words, $f*_\sigma g$ is well-defined element in $L^r(G)$ for $f\in L^p(G)$ and $g\in L^q(G)$ with $\frac{1}{p}+\frac{1}{q} = \frac{1}{r} + 1$, so that we have
				\begin{equation}\label{eq-Young}
				\|f*_\sigma g \|_r \le \|f\|_p \|g\|_q.
				\end{equation}
			
\item Suppose two 2-cocycles $\sigma_1, \sigma_2 \in Z^2(G,\tor)$ are similar,
then we can easily see that the map $f\mapsto f\cdot \xi$ becomes an isometry on $L^p(G)$, $1\le p\le \infty$. Moreover, it is immediate to check that
$$ (f*_{\sigma_1} g)\cdot \xi = (f\cdot \xi)*_{\sigma_2} (g\cdot \xi),\;\; (f\cdot \xi)^{\star_2} = f^{\star_1}\cdot \xi,\;\; f,g\in C_c(G),$$
where $\star_1$ and $\star_2$ are the twisted involutions corresponding to $\sigma_1$ and $\sigma_2$, respectively.
Thus, we can see that the map $L^1(G,\sigma_1) \to L^1(G,\sigma_2)$, $f\mapsto f\cdot \xi$ becomes an isometric $*$-isomorphism.
\end{enumerate}	
\end{rem}

The ``twisting" structure on $L^1(G,\sigma)$ can be transferred to $L^1(G_\sigma)$, an ordinary convolution algebra, via the 
following embedding and projection.
$$j: L^1(G,\sigma) \to L^1(G_\sigma),\;\; f\mapsto j(f),\;\; P: L^1(G_\sigma) \to L^1(G, \sigma),\;\; F \mapsto PF,$$
where $j(f)(s,z) := \bar{z}f(s)$, $(s,z) \in G_\sigma$ and $PF(s) := \int_\tor F(s,z)zdz$, $s\in G$.
Note that it is straightforward to see that $j$ is an isometry and $P$ is a contraction, and both of $j$ and $P$ are $*$-homomorphisms such that $P(j(f)) = f$, $f\in L^1(G,\sigma)$.
This allows us view the twisted convolution algebra $L^1(G,\sigma)$ as a part of the usual convolution algebra $L^1(G_\sigma)$.
This correspondence allows us a well-behaving bounded approximate identity on $L^1(G_\sigma)$, which we record in the following proposition (see also \cite[Lemma~3.1]{Edw}).

\begin{prop}\label{prop-BAI}
Let $\sigma \in Z^2(G,\tor)$ be locally continuous at the identity,
then there is a bounded approximate identity
$(\psi_U)_{U \in \mathcal{U}}$ of $L^1(G,\sigma)$ such that ${\rm supp}\,\psi_U \subseteq U$ and $\|\psi_U\|_{L^1(G,\sigma)}\le 1$ for each $U \in \mathcal{U}$. 
\end{prop}
\begin{proof}
Let $V$ be a neighborhood of the identity where $\sigma$ is continuous on $V\times V$
and $\mathcal{U}$ be a neighborhood base at the identity of $G$ included in $V$.
Then, we select $n_U \in \n$ such that $n_U \to \infty$ as $U \to \{e\}$, where $e$ is the identity of $G$.
Since $\sigma$ is continuous on $V\times V$ we know that $(U \times (-\frac{1}{n_U}, \frac{1}{n_U}))_{U \in \mathcal{U}}$ is a neighborhood base at the identity of $G_\sigma$ in the Weil topology.
Here, we identified $\tor \cong [-\pi, \pi]$.
By the usual construction of bounded approximate identity for $L^1(G_\sigma)$ we get a family of functions $\Psi_U$ such that ${\rm supp}\Psi_U \subseteq U \times (-\frac{1}{n_U}, \frac{1}{n_U})$ and $\|\Psi_U\|\le 1$ for each $U \in \mathcal{U}$.
Now we set
\begin{equation}\label{eq-BAI}
\psi_U := P(\Psi_U),\;\; U \in \mathcal{U}.
\end{equation}
Then, it is straightforward to check that $(\psi_U)_{U \in \mathcal{U}}$ is a bounded approximate identity of $L^1(G,\sigma)$ we wanted.
\end{proof}

We can define the {\it twisted full group $C^*$-algebra} $C^*(G,\sigma)$ by the enveloping $C^*$-algebra of $L^1(G,\sigma)$, but we need more details about twisted representations, which will be essential for the next steps.

\begin{defn}\label{def-sigma-rep}
Let $\sigma \in Z^2(G,\tor)$. A map $\pi : G \to \mathcal{U}(\Hi)$ for some Hilbert space $\Hi$ is called a {\bf $\sigma$-projective unitary representation} (shortly, {\bf $\sigma$-representation}) on $G$ if the function $s \in G \mapsto \la \pi(s)\xi, \eta \ra$ is Borel measurable for any $\xi, \eta\in \Hi$ and we have
\begin{equation}\label{eq-proj-rep}
\pi(s)\pi(t) = \sigma(s,t)\pi(st)
\end{equation}
for any $s,t\in G$.
Here, $\mathcal{U}(\Hi)$ refers to the group of unitary operators on $\Hi$,
and the above function $\pi_{\xi,\eta}(s) := \la \pi(s)\xi, \eta \ra,\; s\in G$ is called the {\bf coefficient function}.
\end{defn}

\begin{rem}\label{rem-tp-rep}
Let $\sigma_i \in Z^2(G,\tor)$ and $\pi_i : G \to B(\Hi_i)$ be $\sigma_i$-representation for $i=1,2$. Then it is straightforward to see that $\sigma_1 \sigma_2 \in Z^2(G,\tor)$ and the tensor product $\pi_1 \otimes \pi_2$ is a $\sigma_1 \sigma_2$-representation.
\end{rem}

One distinctive example of $\sigma$-representation is the left regular one as follows.
\begin{defn}
Let $\sigma \in Z^2(G,\tor)$. The {\bf left regular $\sigma$-representation} $\lambda_\sigma$  is given by $\lambda_\sigma: G\to \mathcal{U}(L^2(G))$ and
$$\lambda_\sigma(s)f(t) := \sigma(s,s^{-1}t)f(s^{-1}t),\; \forall s,t \in G,\; f\in L^2(G).$$
The above operator $\lambda_\sigma(s)$, $s\in G$ can also be understood as an isometry on $L^p(G)$, $1\le p \le \infty$.
\end{defn}

\begin{rem}\label{rem-left-regular-rep-twisted-conv}
\begin{enumerate}
\item Suppose $f,g\in L^2(G)$ with $f$ being compactly supported, in other words, there is a compact set $K\subseteq G$ such that $f = 0$ almost everywhere on $K^c$.
Then, we have
$$\la \lambda_\sigma(s)f, g \ra = \overline{g *_\sigma \tilde{f}(s)}, \;s\in G,$$
where
\begin{equation}\label{eq-tilde}
\tilde{f}(t) := \overline{\sigma(t,t^{-1})}\overline{f(t^{-1})} = \Delta_G(t)f^\star(t),\;\; t\in G.
\end{equation}
This explains that the function $s\in G\mapsto \la \lambda_\sigma(s)f, g\ra$ is Borel measurable.
Passing through limits we can conclude that the function $s\in G\mapsto \la \lambda_\sigma(s)f, g\ra$ is Borel measurable for any $f, g\in L^2(G)$,
which means that $\lambda_\sigma$ is indeed a $\sigma$-representation.
\item It is straightforward to check that $\lambda_\sigma(s)$, $s\in G$, commutes with the left twisted convolution operation on $L^1(G)$, i.e.
\begin{equation}\label{eq-commuting}
(\lambda_\sigma(s)f)*_\sigma g = \lambda_\sigma(s)(f*_\sigma g),\;\; f,g\in L^1(G).
\end{equation}
\end{enumerate}
\end{rem}

As in the untwisted case any $\sigma$-representation $\pi: G\to \mathcal{U}(\Hi)$ can be lifted to a $*$-representation $\tilde{\pi}: L^1(G,\sigma) \to B(\Hi)$ given by
$$\tilde{\pi}(f) := \int_G f(s)\pi(s) ds.$$
For the converse direction we need the non-degeneracy of $\tilde{\pi}$.
Recall that $\tilde{\pi}$ is called {\it non-degenerate} if the subspace $\text{span}\,\tilde{\pi}(L^1(G)) \Hi $ is dense in $\Hi$,
or equivalently for each non-zero $\xi \in \Hi$ there is $f\in L^1(G)$ such that $\tilde{\pi}(f)\xi \ne 0$.
This requires a local continuity of $\pi$, which can be obtained by an interplay between $(G,\sigma)$ and $G_\sigma$.

\begin{prop}\label{prop-loc-continuity}
Let $\sigma \in Z^2(G,\tor)$ be locally continuous at the identity.
For a $\sigma$-representation $\pi: G\to \mathcal{U}(\Hi)$ acting on a separable Hilbert space $\Hi$ there is a neighborhood $V$ of the identity such that
the function $s \in G \mapsto \la \pi(s)\xi, \eta \ra$ is continuous on $V$ for any $\xi, \eta \in \Hi$.
\end{prop}
\begin{proof}
We follow the argument in \cite[p.124]{EL}.
Let $\Pi: G_\sigma \to \mathcal{U}(\Hi)$ be the associated representation of $G_\sigma$ given by 
$$\Pi(s,z) := z\pi(s),\;\;(s,z)\in G_\sigma.$$
It is a homomorphism, which is weakly Borel measurable. 
By \cite[(22.20)]{HR1} we know that $\Pi$ is indeed a strongly continuous $*$-representation since $\Hi$ is separable.
Now we only need to recall \eqref{eq-local-conti-embed} for the conclusion.
\end{proof}

In particular, since the group $G$ is assumed to be second countable, the left regular $\sigma$-representation $\lambda_\sigma$ is weakly locally continuous at the identity, meaning that the function $s \in G \mapsto \la \lambda_\sigma(s)\xi, \eta \ra$ is continuous on a neighborhood $V$ of the identity for any $\xi, \eta \in  L^2(G)$.

\begin{prop}\label{prop-lifting-rep}
Let $\sigma \in Z^2(G,\tor)$ be locally continuous at the identity.
\begin{enumerate}
\item For a $\sigma$-representation $\pi: G\to \mathcal{U}(\Hi)$ on $G$ the lifted map $\tilde{\pi}: L^1(G,\sigma) \to B(\Hi)$
is a $*$-representation satisfying
\begin{equation}\label{eq-rep-lift}
\pi(s)\tilde{\pi}(f) = \tilde{\pi}(\lambda_\sigma(s)f),\; s\in G.
\end{equation}
Moreover, if $\Hi$ is separable, then $\tilde{\pi}$ is non-degenerate.
			
\item For a non-degenerate $*$-representation $\Phi: L^1(G,\sigma) \to B(\Hi)$ we can find a uniquely determined $\sigma$-representation $\pi: G\to \mathcal{U}(\Hi)$ such that $\Phi = \tilde{\pi}$.
\end{enumerate}	
\end{prop}
\begin{proof}
(1) Checking that $\tilde{\pi}$ is a $*$-representation and \eqref{eq-rep-lift} comes from a straightforward calculation.  

When $\Hi$ is separable, we can follow a standard argument for the non-degeneracy of $\tilde{\pi}$ together with the local continuity.
Let $V$ be a neighborhood of the identity such that the function $s \in G \mapsto \la \pi(s)\xi, \eta \ra$ is continuous on $V$ for any $\xi, \eta \in \Hi$.
For any nonzero $\xi\in \Hi$ we take a neighborhood $U\subseteq V$ of the identity such that
$$\|\pi(s)\xi - \xi\| < \frac{\|\xi\|}{2},\;\;s\in U.$$
Note that in the above we are using
$$\|\pi(s)\xi - \xi\|^2 = \|\pi(s)\xi\|^2 + \|\xi\|^2 - \la \pi(s)\xi, \xi\ra - \la \xi, \pi(s)\xi\ra$$
and $\|\pi(s)\xi\| \le \|\xi\|$.
Now we set $f = \frac{1}{\mu(U)}1_U$, where $\mu = ds$ is the fixed left Haar measure on $G$.
Then, we have
$$\|\tilde{\pi}(f)\xi - \xi\|\le \frac{1}{\mu(U)}\int_U \|\pi(s)\xi - \xi\| ds \le \frac{\|\xi\|}{2} <\|\xi\|,$$
which means that $\tilde{\pi}(f)\xi\ne 0$.

(2) We follow the same argument as in the untwisted case. Let $(\psi_U)_{U\in \mathcal{U}}$ be the bounded approximate identity of $L^1(G,\sigma)$ from Proposition \ref{prop-BAI}, then for any $s\in G$ and $f\in L^1(G)$ we have
$$(\lambda_\sigma(s)\psi_U)*_\sigma f = \lambda_\sigma(s)(\psi_U*_\sigma f) \to \lambda_\sigma(s)f\;\; \text{in}\;\; L^1(G)$$
by \eqref{eq-commuting}, which implies
$$\Phi(\lambda_\sigma(s)\psi_U)\Phi(f)\xi \to \Phi(\lambda_\sigma(s)f)\xi,\;\; \xi\in \Hi.$$
This allows us to define the operator $\pi(s)$ given by
\begin{equation}\label{eq-rep-def}
\pi(s)\Phi(f)\xi := \Phi(\lambda_\sigma(s)f)\xi,\;\; \xi\in \Hi,
\end{equation}
which is clearly a contraction on ${\rm span}\, \Phi(L^1(G))\Hi$. Since it is a dense subspace of $\Hi$, thanks to the non-degeneracy of $\Phi$, the map $\pi(s)$ extends to a contraction on $\Hi$.
Now, it is straightforward to check that $\pi$ is a $*$-homomorphism, which, in turn, implies that $\pi(s)$ is a unitary for all $s\in G$.

Finally, we focus on the measurability.
The same argument as in Remark \ref{rem-left-regular-rep-twisted-conv} (1) tells us that that the function $s\in G \mapsto \la\lambda_\sigma(s)f, g\ra$ is Borel measurable for any $f\in L^1(G)$ and $g\in L^\infty(G)$.
Thanks to the separability of $L^2(G)$ we can appeal to the Pettis measurability theorem (\cite[Theorem~II.2]{DU}) to see that the map $s\in G\mapsto \lambda_\sigma(s)f \in L^1(G)$ is Bochner measurable.
Since $\Phi$ is continuous, we know that the map 
$$s\in G\mapsto \Phi(\lambda_\sigma(s)f) = \pi(s)\Phi(f) \in B(\Hi)$$
is Bochner measurable, which means that the map $s\in G\mapsto \la \pi(s)\Phi(f)\xi, \eta\ra$ is Borel measurable for any $\xi, \eta \in \Hi$.
We use the density of ${\rm span}\, \Phi(L^1(G))\Hi$ in $\Hi$ and the limit procedure to get the conclusion we wanted.
\end{proof}

\subsection{Twisted group $C^*$-algebras and von Neumann algebras}

We begin with the universal group $C^*$-algebra with twisting.
\begin{defn}
Let $\sigma \in Z^2(G,\tor)$ be locally continuous at the identity. We define the {\bf full twisted group $C^*$-algebra $C^*(G,\sigma)$} of $G$ to be the enveloping $C^*$-algebra of the Banach $*$-algebra $L^1(G,\sigma)$. In other words, for any $f\in L^1(G)$ we have the associated norm
\begin{equation}\label{eq-full-norm}
\|f\|_{C^*(G,\sigma)} := \sup\{\|\tilde{\pi}(f)\|_{B(\Hi)}\;|\; \pi: G\to B(\Hi),\; \text{$\sigma$-representation}\},
\end{equation}
where $\tilde{\pi}$ is the lifted representation as in Proposition \ref{prop-lifting-rep} and $C^*(G,\sigma)$ is the completion of $(L^1(G), \| \cdot \|_{C^*(G,\sigma)})$. The canonical embedding from $L^1(G,\sigma)$ into $C^*(G,\sigma)$ will be denoted by $\iota$, i.e.
\begin{equation}\label{eq-full-embed}
\iota : L^1(G,\sigma) \hookrightarrow C^*(G,\sigma).
\end{equation}
\end{defn}

\begin{rem}
The above norm \eqref{eq-full-norm} is actually a well-defined $C^*$-norm thanks to the fact that the lifted left regular $\sigma$-representation $\widetilde{\lambda_\sigma}$ is injective. Indeed, for any $f\in L^1(G)$ with $\widetilde{\lambda_\sigma}(f) = 0$ we have
$$f = \lim_{U} f*_\sigma \psi_U = \lim_U \widetilde{\lambda_\sigma}(f)\psi_U = 0,$$
where $(\psi_U)_{U \in \mathcal{U}}$ is the bounded approximate identity of $L^1(G,\sigma)$ in Proposition \ref{prop-BAI}, so that we know $\psi_U \in L^2(G)$ for each $U$.
\end{rem}

Now we move to the reduced version and its von Neumann algebra counterpart.
	\begin{defn}
		Let $\sigma \in Z^2(G,\tor)$. We define the {\bf reduced twisted group $C^*$-algebra $C^*_r(G,\sigma)$} of $G$ to be the $C^*$-algebra generated by $\{\widetilde{\lambda_\sigma}(f) : f\in L^1(G)\}$ in $B(L^2(G))$. We also define the {\bf twisted group von Neumann algebra $VN(G,\sigma)$} of $G$ to be the von Neumann algebra generated by $\{\widetilde{\lambda_\sigma}(f) : f\in L^1(G)\}$ (or equivalently by $\{\lambda_\sigma(s) : s\in G\}$) in $B(L^2(G))$.
	\end{defn}
	
	\begin{ex}\label{ex-Q-torus-plane}
		\begin{enumerate}
			\item Let $A = \begin{bmatrix}0 & \theta\\-\theta & 0\end{bmatrix}$ for an irrational $\theta \in [0,1]$. Then, for the 2-cocycle $\sigma_A$ from \eqref{eq-2-cocyles-ex} we have that $C^*(\z^2, \sigma_A)$ is the non-commutative torus $C(\tor^2_\theta)$ as mentioned in the introduction.
			
			\item Let $B = \begin{bmatrix}0 & I_n\\-I_n & 0\end{bmatrix} \in M_{2n}(\Real)$, where $I_n$ is the $n\times n$ identity matrix. Then, $VN(\Real^{2n}, \sigma_B) \cong B(L^2(\Real^n))$ as von Neumann algebras \cite[example p.490]{KL72}.
		\end{enumerate}
	\end{ex}
	
\section{Twisted Fourier(-Stieltjes) spaces}\label{sec-twisted-space}

In this section we define twisted Fourier(-Stieltjes) spaces $A(G,\sigma)$ and $B(G,\sigma)$ extending the discrete group case of B\'{e}dos/Conti \cite{BC2009, BC2016-2}.
Note that we assume the local continuity at the identity of the associated 2-cocycle from now on.
We mainly follow the classical approach to define Fourier(-Stieltjes) algebras in the presence of twisting. We begin with twisted ``positive definite" functions.

\begin{defn}
Let $\sigma \in Z^2(G,\tor)$  be locally continuous at the identity. 
We say that $\varphi: G \to \Comp$ is of {\bf $\sigma$-positive type} if $\varphi \in L^\infty(G)$ and it defines a positive linear functional on the Banach $*$-algebra $L^1(G,\sigma)$, i.e. for any $f\in L^1(G)$ we have
	$$\int_G f^\star *_\sigma f(s) \varphi(s) ds \ge 0.$$
\end{defn}

\begin{rem}
We may define $\varphi: G \to \Comp$ to be a {\bf $\sigma$-positive definite} function if for any finitely supported function $\alpha: G\to \Comp$ we have
$$\sum_{s,t\in G}\overline{\sigma(s,s^{-1}t)}\varphi(s^{-1}t)\overline{\alpha(s)}\alpha(t) \ge 0.$$
However, without continuity of $\sigma$ we actually need to work on functions $\varphi$ of $\sigma$-positive type only with local continuity at the identity.
For those $\varphi$ it is not clear whether they are $\sigma$-positive definite.

When $\sigma$ is continuous, for a bounded continuous function $\varphi: G \to \Comp$, it is of $\sigma$-positive type if and only if it is $\sigma$-positive definite.
Here comes the proof of this. Note that we can readily check
$$\int_G f^\star *_\sigma f(s) \varphi(s) ds = \int_G\int_G \overline{\sigma(s,s^{-1}t)}\varphi(s^{-1}t)\overline{f(s)}f(t) dsdt,$$
so that we may use the same argument as in the untwisted case.
One key observation is that for any compact $K \subseteq G$ we have $M(K) = C(K)^*$ is generated by the span of point masses on $K$ or $L^1(K)$ in the weak$^*$ topology.
\end{rem}

\begin{defn}
Let $\sigma \in Z^2(G,\tor)$ be locally continuous at the identity.
We define $\mathcal{P}(G,\sigma)$ to be the collection of functions of $\sigma$-positive type, which is locally continuous at the identity. 
We also define
$$\mathcal{P}_1(G,\sigma) := \{\varphi \in \mathcal{P}(G,\sigma): \varphi(e)=1\}.$$
\end{defn}

\begin{rem}\label{rem-ptwise-eval}
The expression ``$\varphi \in \mathcal{P}(G,\sigma)$" means that we can choose a representative function of the equivalence class (upto almost everywhere coincidence), which is locally continuous at the identity.
Note that the choice of such representative function is locally unique around the identity.
If we denote the chosen representative function by $\varphi$ again by abuse of notation,
then the pointwise evaluation $\varphi(e)$ is well defined.
Moreover, for $f\in L^1(G)$ we have
$$\int_G f^\star *_\sigma f(s) \varphi(s) ds = \int_G\int_G \overline{\sigma(s,s^{-1}t)}\varphi(s^{-1}t)\overline{f(s)}f(t) dsdt.$$
Setting $f = \frac{1}{|V|}1_V$ for a neighborhood $V$ of the identity with the Haar measure $|V|$ being finite and taking $V \to \{e\}$ we get $\varphi(e)\ge 0$.
\end{rem}

Now we present a characterization of $\mathcal{P}(G,\sigma)$ using coefficient functions from Definition \ref{def-sigma-rep}.

\begin{prop}\label{prop-P1-equiv}
Let $\sigma \in Z^2(G,\tor)$ be locally continuous at the identity and $\varphi\in L^\infty(G)$.
Then, the following are equivalent.
\begin{enumerate}
\item $\varphi \in \mathcal{P}_1(G,\sigma)$.
\item The functional $f \in L^1(G)\mapsto \int_G f(s)\varphi(s)ds$ determines a state $\tilde{\varphi}$ on $C^*(G,\sigma)$.
\item There is a $\sigma$-representation $\pi: G \to B(\Hi)$ acting on a separable Hilbert space $\Hi$ \color{black}and a norm 1 element $\xi\in \Hi$ such that $\varphi(s) = \la \pi(s)\xi, \xi \ra$ for almost every $s\in G$.\color{black}
\end{enumerate}
\end{prop}
\begin{proof}
(1) $\Rightarrow$ (2) The positivity of the associated functional on the Banach $*$-algebra $L^1(G,\sigma)$ is from the assumption on $\varphi$.
Now we recall the bounded approximate identities $(\psi_U)_{U \in \mathcal{U}}$ and $(\Psi_U)_{U \in \mathcal{U}}$ of $L^1(G,\sigma)$ and $L^1(G_\sigma)$, respectively,
in the proof of Proposition \ref{prop-BAI}. 
Then, we have
\begin{align*}
\lim_U\int_G \psi_U(s)\varphi(s)ds
& = \lim_U\int_G \int_\tor \Psi_U(s,z)\varphi(s)zdzds\\
& = \lim_U\int_G \int_\tor \Psi_U(s,t)J(\varphi)(s,z)dzds\\
& = J(\varphi)(e,1) = \varphi(e) = 1,
\end{align*}
where $J(\varphi)(s,z) = \varphi(s)z$, $(s,z) \in G\times \tor$.
This explains that the associated functional extends to a state on $C^*(G,\sigma)$, which is the enveloping $C^*$-algebra of $L^1(G,\sigma)$.

\vspace{0.3cm}

(2) $\Rightarrow$ (3) Let $(\Phi, \Hi, \xi)$ be the GNS representation of $C^*(G,\sigma)$ with respect to the state $\tilde{\varphi}$.
Note that $\Hi$ is separable since $G$ is second countable.
By restricting to $L^1(G,\sigma)$ we get a non-degenerate $*$-representation of $L^1(G,\sigma)$, and in turn by Proposition \ref{prop-lifting-rep} (2) we get a $\sigma$-representation $\pi: G \to B(\Hi)$ such that 
$\Phi \circ \iota = \tilde{\pi}$, where $\iota: L^1(G,\sigma) \hookrightarrow C^*(G,\sigma)$ is the canonical embedding and $\tilde{\pi}$ is the lifting of $\pi$. Thus, we have
	$$\tilde{\varphi}(\iota(f)) = \la \Phi(\iota(f))\xi, \xi \ra = \la \tilde{\pi}(f)\xi, \xi \ra =  \int_G f(s) \la \pi(s)\xi, \xi \ra ds,$$
so that we get $\varphi(s) = \la \pi(s)\xi, \xi \ra$, $s\in G$.	

\vspace{0.3cm}

(3) $\Rightarrow$ (1) If $\varphi(s) =\langle \pi (s) \xi, \xi \rangle $ for some $\sigma$-representation $\pi: G\to B(\Hi)$ on a separable Hilbert space $\Hi$ and $\xi \in \Hi$ with norm 1,
then $\varphi \in L^\infty(G)$ is locally continuous at the identity by Proposition \ref{prop-loc-continuity} and clearly $\varphi(e) = 1$. It remains to check that $\varphi$ is of $\sigma$-positive type. By change of variables we have for any $f\in C_c(G)$ that
\begin{align*}
\int_G f^\star  *_\sigma f(s) \varphi(s) ds &= \int _G \int_G \Delta_G(t^{-1})  \overline{\sigma(t,t^{-1})}  \sigma(t, t^{-1}s) \overline{f(t^{-1})}   f(t^{-1}s )  \varphi(s) \, ds\, dt\\
&= \int _G \int_G \Delta_G(t^{-1})  \overline{\sigma(t,t^{-1})}  \sigma(t,  s) \overline{f(t^{-1})}   f(s )  \varphi(ts) \, ds\, dt\\
&= \int _G \int_G \Delta_G(t^{-1})  \overline{\sigma(t,t^{-1})}  \,\overline{f(t^{-1})}   f(s )  \langle \pi (t) \pi (s) \xi, \xi \rangle \, ds\, dt\\
&= \int _G \int_G   \overline{\sigma(t,t^{-1})}  \,\overline{f(t)}   f(s )  \langle \pi (t^{-1}) \pi (s) \xi, \xi \rangle \, ds\, dt\\
&= \int _G \int_G   \overline{f(t)}   f(s )  \langle  \pi (s) \xi, \pi (t)\xi \rangle \, ds\, dt\\
&=   \langle  \pi (f) \xi, \pi (f)\xi \rangle  \geq 0.
\end{align*}
We may use density for general $f\in L^1(G)$.
\qedhere
\end{proof}

\begin{cor}\label{cor-sigma-pd}
Let $\sigma_1, \sigma_2$ and $\sigma \in Z^2(G,\tor)$ be locally continuous at the identity.
\begin{enumerate}
\item If $\varphi \in \mathcal{P}_1(G,\sigma)$, then $\overline{\varphi} \in \mathcal{P}_1(G, \overline{\sigma})$.
		
\item If $\varphi_1 \in \mathcal{P}_1(G,\sigma_1)$ and $\varphi_2 \in \mathcal{P}_1(G,\sigma_2)$, then we have $\varphi_1\varphi_2 \in \mathcal{P}_1(G,\sigma_1 \sigma_2)$.
		
\item For compactly supported $f \in L^2(G)$ we have
$$\overline{f *_\sigma \tilde{f}} \in \mathcal{P}(G,\sigma),$$
where $
\tilde{f}(t) := \overline{\sigma(t,t^{-1})}\overline{f(t^{-1})} = \Delta_G(t)f^\star(t),\; t\in G.$ 
\end{enumerate}
\end{cor}
\begin{proof}
(1) is immediate from the definition of positive definite functions and (2) is
clear from Remark \ref{rem-tp-rep}.
Finally, (3) is directly from \eqref{eq-tilde} and Proposition \ref{prop-P1-equiv}.
\end{proof}

\begin{prop}\label{prop-cpt-supp-PD}
Let $\sigma \in Z^2(G,\tor)$ be locally continuous at the identity  and $\varphi \in \mathcal{P}(G,\sigma)$ with compact support.
Then, there is $\psi \in L^2(G)$ such that $\varphi(s) = \la \lambda_\sigma(s)\psi, \psi \ra$ for almost every $s\in G$.
\end{prop}
\begin{proof}

We follow the argument of \cite[Section 13.8]{Dix} closely. We first observe that, for a compactly supported $\varphi \in L^\infty(G)$, the function $\varphi$ is $\sigma$-positive definite if and only if the associated right twisted convolution operator
$$R(\bar{\varphi}) : L^2(G) \to L^2(G),\; f \mapsto f*_\sigma \bar{\varphi}$$
is positive. Indeed, we have
\begin{align*}
\la f*_\sigma \bar{\varphi}, f \ra
& = \int_G\int_G f(t) \sigma(t,t^{-1}s)\overline{\varphi(t^{-1}s)}\overline{f(s)}dtds\\
& = \overline{\int_G f^\star  *_\sigma f(s) \varphi(s) ds}
\end{align*}
for any $f\in C_c(G) \subseteq  L^2(G)$.
Multiplying by a positive constant we may assume that $0\le R(\bar{\varphi}) \le 1$.

Now we consider an increasing sequence $(p_n(t))_{n\ge 1}$ of non-negative polynomials on $[0,1]$, vanishing at 0 and converging uniformly to the function $\sqrt{t}$ on $[0,1]$. We set
$$\psi_j := p_j(\bar{\varphi}) \in L^1(G,\sigma)$$
using the functional calculus of the Banach $*$-algebra $L^1(G,\sigma)$. Then the function $\psi_j$ is compactly supported as well and $R(\psi_j) = p_j(R(\bar{\varphi}))$, so that
$$0\le R(\psi_1) \le R(\psi_2) \le \cdots \le R(\bar{\varphi})\le 1.$$
Moreover, we can see that $\{R(\psi_j): j \ge 1\} \cup \{R(\bar{\varphi})\}$ is a commuting family of positive operators.
Note that for any compactly supported $g,h\in L^\infty(G)$ we have $R(g)R(h) = R(g *_\sigma h)$ and $R(g)^* = R(\tilde{g})$, where $\tilde{g}$ is from \eqref{eq-tilde}.
Using this observation we have
$$R(\psi_k *_\sigma \tilde{\psi_j}) = R(\psi_k)R(\psi_j)^* = R(\psi_k)R(\psi_j)\ge 0,\;\; \forall j,k\ge 1.$$
By Corollary \ref{cor-sigma-pd} we know that $\overline{\psi_k *_\sigma \tilde{\psi_j}}\in \mathcal{P}(G,\sigma)$ and consequently we have
$$\la \psi_k, \psi_j \ra = \psi_k *_\sigma \tilde{\psi_j}(e) \ge 0,\;\; \forall j,k\ge 1.$$
Note that $\psi_k *_\sigma \tilde{\psi_j} = \overline{\la \lambda_\sigma(s)\psi_j, \psi_k\ra}$, $s\in G$, so that we have chosen a representative function which is locally continuous at the identity by Proposition \ref{prop-loc-continuity}.
Thus, the evaluation $\psi_k *_\sigma \tilde{\psi_j}(e)$ is well-defined and non-negative as in Remark \ref{rem-ptwise-eval}.
Similarly, we can observe that $\la \psi_k - \psi_j, \psi_j  \ra \ge 0$ since
$R((\psi_k - \psi_j) *_\sigma \tilde{\psi_j}) = [R(\psi_k) - R(\psi_j)]R(\psi_j)\ge 0$ for any $k\ge j\ge 1.$ 

Combining all the above for any $k\ge j\ge 1$ we have
\begin{equation}\label{eq-Cauchy}
\|\psi_k - \psi_j\|^2_2 = \|\psi_k\|^2_2 - \|\psi_j\|^2_2 - 2\la \psi_k - \psi_j, \psi_j  \ra \le \|\psi_k\|^2_2 - \|\psi_j\|^2_2,
\end{equation}
so that we have $\|\psi_j\|_2 \le \|\psi_k\|_2$. By replacing $\psi_k$ with $\bar{\varphi}$, the same argument tells us that $\|\psi_j\|_2 \le \|\bar{\varphi}\|_2$, which implies that
	$$\|\psi_1\|_2 \le \|\psi_2\|_2 \le \cdots \le \|\bar{\varphi}\|_2.$$
Thus, we know that the sequence $(\|\psi_j\|_2)_{j\ge 1}$ converges and \eqref{eq-Cauchy} implies that $(\psi_j)_{j\ge 1}$ is a Cauchy sequence in $L^2(G)$. Thus, we get $\displaystyle \psi = \lim_{j\to \infty}\psi_j\in L^2(G)$. Finally, for $h\in C_c(G)$ we have
\begin{align*}
h *_\sigma \bar{\varphi}
& = R(\bar{\varphi})(h) = \lim_{j\to \infty} p^2_j(R(\bar{\varphi}))(h)\\
& =  \lim_{j\to \infty} R(\psi_j *_\sigma \tilde{\psi_j})(h) = \lim_{j\to \infty} h*_\sigma \psi_j *_\sigma \tilde{\psi_j}\\
& = h*_\sigma \psi *_\sigma \tilde{\psi},
\end{align*}
where the last equality is due to Young's inequality \eqref{eq-Young}. This means that $\varphi(s) = \overline{\psi *_\sigma \tilde{\psi}}(s)$ for almost every $s\in G$.
\end{proof}

\begin{defn}
Let $\sigma \in Z^2(G,\tor)$ be locally continuous at the identity. We define the {\bf twisted Fourier-Stieltjes space} $B(G,\sigma)$ by
$$B(G,\sigma) := \{ \pi_{\xi,\eta}\; | \;\pi: G \to B(\Hi)\; \text{$\sigma$-representation, $\Hi$ separable},\; \xi, \eta \in \Hi \},$$
where $\pi_{\xi,\eta}$ is a coefficient function as in Definition \ref{def-sigma-rep}.	
\end{defn}

\begin{rem}
\begin{enumerate}
\item
From Proposition \ref{prop-P1-equiv} and Proposition \ref{prop-loc-continuity} it is immediate to see that $B(G,\sigma) = {\rm span}\,\mathcal{P}_1(G,\sigma)$ and thus we have the isometric identification $B(G,\sigma) = \left( C^*(G,\sigma)\right )^*$ if we equip $B(G,\sigma)$ with the norm
$$\|\varphi\|_{B(G,\sigma)} := \sup_{f\in L^1(G),\, ||\iota(f)||_{C^*(G,\sigma)} \le 1} \left| \int_G f(s)\varphi(s)ds \right|,$$
where $\iota$ is the canonical embedding in \eqref{eq-full-embed},
which gives us the contractive embedding
$$B(G,\sigma) \hookrightarrow L^\infty(G).$$
Note that the norm $\|\varphi\|_{B(G,\sigma)}$ is determined by the values $\int_G f(s)\varphi(s)ds$ for all $f\in L^1(G)$, which justifies the inclusion $B(G,\sigma) \subseteq L^\infty(G)$.
\item We know from Proposition \ref{prop-loc-continuity} that functions in $B(G,\sigma)$ are locally continuous at the identity.  
However, we still need the global continuity of $\sigma$ to ensure that $B(G,\sigma) \subseteq C_b(G)$, the space of bounded continuous functions on $G$. 		
\item In \cite{BC2016} a Fourier-Stieltjes algebra $B(\Sigma)$ for a $C^*$-dynamical system $\Sigma = (A, G, \alpha, \sigma)$ has been introduced.
Here, $A$ is a unital $C^*$-algebra, $G$ is a discrete group, $\alpha: G \to {\rm Aut}(A)$ and $\sigma: G\times G \to \mathcal{U}(A)$ satisfying a certain compatibility condition.
If we take the $A = \Comp$ and $\alpha_g \equiv id$, $g\in G$, then $\sigma$ is a 2-cocycle in our sense. However, for an equivariant representation $(\rho, v)$, we can actually see that $v$ does not depend on $\sigma$, so that it becomes a usual representation, at best. 
This means that the algebra $B(\Sigma)$ does not cover the case of the Fourier-Stieltjes space $B(G,\sigma)$ introduced in this article.
\end{enumerate}
\end{rem}

The dual space $B(G,\sigma)^*$ of $B(G,\sigma)$ is the bidual of $C^*(G,\sigma)$, which can be understood as the universal enveloping von Neumann algebra of $C^*(G,\sigma)$ in a canonical way.
Let us denote the mentioned algebra by $W^*(G,\sigma)$.
More precisely, there is the universal representation
$$\Phi^\sigma_u : C^*(G,\sigma) \to B(\Hi^\sigma_u),$$
which is a direct sum of all cyclic representations such that $W^*(G,\sigma) = \Phi^\sigma_u(C^*(G,\sigma))^{''}$.
Note that we have
\begin{equation}\label{eq-universal-rep-duality}
\la \Phi^\sigma_u(T), \varphi \ra = \la T, \varphi \ra,\;\; T\in C^*(G,\sigma), \;\varphi \in B(G,\sigma).
\end{equation}
By Proposition \ref{prop-lifting-rep} there is a $\sigma$-representation
$\pi^\sigma_u: G\to B(\Hi^\sigma_u)$
called the {\bf twisted universal representation} such that its lifted representation is given by 
\begin{equation}\label{eq-twisted-universal-rep}
\widetilde{\pi^\sigma_u} = \Phi^\sigma_u \circ \iota,
\end{equation}
where $\iota : L^1(G,\sigma) \to C^*(G,\sigma)$ is the canonical embedding in \eqref{eq-full-embed}.
Now, we recall the definition of the usual universal enveloping von Neumann algebra $W^*(G)$ of $C^*(G)$ and the (untwisted) universal representation $\pi_u: G \to B(\Hi)$ with $W^*(G) = \widetilde{\pi_u} (L^1(G))^{''}$.
Tensoring the above two we get another $\sigma$-representation
$$\pi_u\otimes \pi^\sigma_u: G \to W^*(G) \bar{\otimes} W^*(G,\sigma) \subseteq B(\Hi_u) \bar{\otimes} B(\Hi^\sigma_u),$$
which extends uniquely to a normal $*$-homomorphism
$$\Gamma^\sigma_u : W^*(G,\sigma) \to W^*(G) \bar{\otimes} W^*(G,\sigma),\;\; \pi^\sigma_u(s) \mapsto \pi_u(s) \otimes \pi^\sigma_u(s),\;\; s\in G$$
by universality. This leads us to the following.

\begin{prop}\label{prop-duality-bimodule}
Let $\sigma \in Z^2(G,\tor)$ be locally continuous at the identity.
\begin{enumerate}
\item For any $\varphi \in B(G,\sigma)$, the two functions $\la \pi^\sigma_u(\cdot), \varphi \ra$ and $\varphi$ agree in $L^\infty(G)$.
\item The following map is a complete contraction.
\begin{equation}\label{eq-multiplier}
m_\sigma : B(G)\prt B(G,\sigma) \to B(G,\sigma),\;\; f\otimes g \mapsto fg.
\end{equation}
Moreover, we have $m^*_\sigma = \Gamma^\sigma_u$.
\end{enumerate}

\end{prop}

\begin{proof}
(1) We recall the bounded approximate identity $(\psi_U)_{U\in \mathcal{U}}$ of $L^1(G,\sigma)$ in Proposition \ref{prop-BAI}
and observe that $\widetilde{\pi^\sigma_u}(\psi_U)$ converges to $1_{W^*(G,\sigma)}$ in the weak$^*$-topology by taking a suitable subnet if necessary.
Since the the weak$^*$-topology and the $\sigma$-weak operator topology on $W^*(G,\sigma)$ coincide we can actually show that $\widetilde{\pi^\sigma_u}(\lambda_\sigma(s)\psi_U)$ converges to $\pi^\sigma_u(s)$, $s\in G$ in the weak$^*$-topology.
Indeed, for any $(\xi_i, \eta_i)_{i\ge 1} \subseteq \Hi^\sigma_u$ such that $\sum_i \| \xi_i \| \cdot \| \eta_i \|<\infty$ we have
\begin{align*}
\sum_i \la \widetilde{\pi^\sigma_u}(\lambda_\sigma(s)\psi_U)\xi_i, \eta_i \ra
& = \sum_i \la \pi^\sigma_u(s)\widetilde{\pi^\sigma_u}(\psi_U)\xi_i, \eta_i \ra\\
& = \sum_i \la \widetilde{\pi^\sigma_u}(\psi_U)\xi_i, \pi^\sigma_u(s)^*\eta_i \ra\\
& \to \sum_i \la \xi_i, \pi^\sigma_u(s)^*\eta_i \ra\\
& = \sum_i \la \pi^\sigma_u(s)\xi_i, \eta_i \ra
\end{align*}
as $U \to \{e\}$.
This implies for any $\varphi \in B(G,\sigma)$ that
$$\la \pi^\sigma_u(s), \varphi \ra = \lim_U \la \widetilde{\pi^\sigma_u}(\lambda_\sigma(s)\psi_U), \varphi \ra,$$
and consequently for $f\in L^1(G)$ we have
\begin{align*}
\int_G f(s)\la \pi^\sigma_u(s), \varphi \ra ds 
& = \lim_U \int_G f(s)\la \widetilde{\pi^\sigma_u}(\lambda_\sigma(s)\psi_U), \varphi \ra ds\\
& = \lim_U \la \widetilde{\pi^\sigma_u}(\widetilde{\lambda_\sigma}(f)\psi_U), \varphi \ra\\
& = \la \widetilde{\pi^\sigma_u}(f), \varphi \ra\\
& = \int_G f(s)\varphi(s) ds.
\end{align*}
Thus, we know that $\la \pi^\sigma_u(\cdot), \varphi \ra$ and $\varphi$ agree in the space $L^\infty(G)$.

(2) This is immediate from the above result and the usual duality $(W^*(G), B(G))$ given by
$$\la \widetilde{\pi_u}(g), \psi \ra = \int_G g(s)\psi(s) ds,\;\; g\in L^1(G),\; \psi \in B(G).$$
\end{proof}	

\begin{rem}
\begin{enumerate}
\item The map $m_\sigma$ in \eqref{eq-multiplier} gives an operator $B(G)$-bimodule structure on $B(G,\sigma)$, which can be regarded as the canonical one.

\item  
In general, for $\varphi \in B(G,\sigma)$, we do not have $\la \pi^\sigma_u(\cdot), \varphi \ra=\varphi(\cdot)$ pointwise. Alternatively, using local continuity, we can find a neighborhood $V$ of the identity such that
$$\la \pi^\sigma_u(s), \varphi \ra = \varphi(s),\;\; s\in V.$$
Indeed, we have
\begin{align*}
\la \pi^\sigma_u(s), \varphi \ra
& = \lim_U \la \widetilde{\pi^\sigma_u}(\lambda_\sigma(s)\psi_U), \varphi \ra\\
& = \lim_U \la \iota(\lambda_\sigma(s)\psi_U), \varphi \ra\\
& = \lim_U \int_G \lambda_\sigma(s)\psi_U(t) \varphi(t) dt\\
& = \lim_U \int_G \sigma(s, s^{-1}t)\psi_U(s^{-1}t) \varphi(t) dt\\
& = \varphi(s)
\end{align*}
for $s\in G$ near the identity
by the local continuity of $\sigma$ and $\varphi$ at the identity.
\end{enumerate}

\end{rem}

%
%


\begin{prop}
Let $\sigma \in Z^2(G,\tor)$ be locally continuous at the identity. For every $\varphi \in B(G,\sigma)$ we have
$$\|\varphi \|_{B(G,\sigma)} = \inf \{\|\xi \| \cdot \|\eta\| \},$$
where the above infimum runs over all possible $\sigma$-representation $\pi: G \to B(\Hi)$ and $\xi, \eta \in \Hi$ such that $\varphi(\cdot) = \la \pi(\cdot)\xi, \eta \ra$. Moreover, the infimum can actually be obtained.
\end{prop}
\begin{proof}
Let $\varphi(\cdot) = \la \pi(\cdot)\xi, \eta \ra$ for some $\sigma$-representation $\pi: G \to B(\Hi)$ and $\xi, \eta \in \Hi$. Then for any $f\in L^1(G)$ we have $|\int_G f(s)\varphi(s)ds| = |\la \tilde{\pi}(f) \xi, \eta \ra| \le \|\tilde{\pi}(f)\| \cdot \|\xi \| \cdot \|\eta\|$, which implies that $\|\varphi \|_{B(G,\sigma)} \le \|\xi \| \cdot \|\eta\|$.

Now we will show that the infimum can actually be obtained. We can view $\varphi$ as a normal functional on the von Neumann algebra $W^*(G,\sigma)$. Thus, we have the polar decomposition $\varphi = u |\varphi|$, where $u\in W^*(G,\sigma)$ is a partial isometry and $|\varphi|$ is a positive element of $B(G,\sigma)$,
which means that there are $\sigma$-representation $\pi: G\to B(\Hi)$ and $\eta \in \Hi$ such that $|\varphi|(\cdot) = \la \pi(\cdot) \eta, \eta \ra$ and $\| |\varphi| \|_{B(G,\sigma)} = \|\eta\|^2$.
The universality of $\pi^\sigma_u$ provides a normal $*$-homomorphism $\rho : W^*(G,\sigma) \to \tilde{\pi}(L^1(G))^{''}$ such that $\rho(\pi^\sigma_u(s)) = \pi(s)$. Then, (1) of Proposition \ref{prop-duality-bimodule} ensures that we have
$$\la T, |\varphi| \ra = \la \rho(T), |\varphi| \ra,\;\; T\in W^*(G,\sigma).$$
From (1) of Proposition \ref{prop-duality-bimodule} again we know for any $s\in G$
\begin{align*}
\varphi(s)
& = \la \pi_u^\sigma(s), \varphi \ra = \la \pi_u^\sigma(s), u |\varphi| \ra\\
& =  \la \pi_u^\sigma(s)u, |\varphi| \ra = \la \rho(\pi_u^\sigma(s)u), |\varphi| \ra\\
& =\la \pi(s) \rho(u), |\varphi| \ra = \la \pi(s) \rho(u)\eta, \eta \ra.
\end{align*}
Now we set $\xi = \rho(u)\eta$, then we have $\varphi(\cdot) = \la \pi(\cdot)\xi, \eta \ra$ and
$$\|\xi\| \le \| \rho(u)\| \cdot \|\eta\| \le \|\eta\|.$$
On the other hand, from the first part of this proof we know that
	$$\|\eta\|^2 = \| |\varphi| \|  = \|\varphi\| \le \|\xi\| \cdot \|\eta\|.$$
Thus, we have $\|\xi \| = \|\eta\|$ and consequently $\|\varphi\| = \|\xi\| \cdot \|\eta\|$.
\end{proof}

Now we define the twisted Fourier space $A(G,\sigma)$.
\begin{defn}
Let $\sigma \in Z^2(G,\tor)$ be locally continuous at the identity. We define the {\bf twisted Fourier space $A(G,\sigma)$} to be the predual $VN(G,\sigma)_*$ of the twisted von Neumann algebra $VN(G,\sigma)$.
\end{defn}

Let us fix $\sigma$ for the rest of this section.
Recall that the embedding $VN(G,\sigma) \subseteq B(L^2(G))$ gives us the canonical complete quotient map $B(L^2(G))_* \to VN(G,\sigma)_*$, which is nothing but the restriction map, so that we have
\begin{align*}
VN(G,\sigma)_*
& = \{ \psi|_{VN(G,\sigma)} : \psi \in B(L^2(G))_* \}\\
& = \{ \psi|_{VN(G,\sigma)} : \psi = \sum_{n\ge 1} \om_{\xi_n, \eta_n},\; \xi_n, \eta_n \in L^2(G),\; \sum_{n\ge 1} \|\xi_n\| \cdot \|\eta_n\| <\infty\}\\
& = \overline{\text{span}}\{(\om_{\xi,\eta})|_{VN(G,\sigma)}: \xi, \eta \in L^2(G)\}
\end{align*}
with the norm
$$\|\varphi\|_{VN(G,\sigma)_*} = \inf \left\{\sum_{n\ge 1} \|\xi_n\| \cdot \|\eta_n\| : \varphi = (\sum_{n\ge 1} \om_{\xi_n, \eta_n})|_{VN(G,\sigma)}\right\}.$$
Here, $\om_{\xi, \eta}$ refers to the normal functional on $B(L^2(G))$ given by
$$\om_{\xi, \eta}(T) = \la T\xi, \eta \ra,\;\; T\in B(L^2(G)).$$	
Since $\{\widetilde{\lambda_\sigma}(f): f\in L^1(G)\}$ is weak$^*$-dense in $VN(G,\sigma)$ any element $\varphi \in VN(G,\sigma)_*$ is determined by its values on $\{\widetilde{\lambda_\sigma}(f): f\in L^1(G)\}$.
Note that for $\xi, \eta\in L^2(G)$ we have
$$\la \widetilde{\lambda_\sigma}(f), \om_{\xi,\eta}\ra 
= \la \widetilde{\lambda_\sigma}(f)\xi,\eta\ra = \int_G f(s)\la \lambda_\sigma(s)\xi,\eta\ra ds$$
so that we get an injective linear map
$$J: A(G,\sigma) \to L^\infty(G),\;\; \om_{\xi, \eta} \mapsto \la \lambda_\sigma(\cdot)\xi, \eta \ra.$$
It is clear that $J$ is a complete contraction and $J(\varphi)$ is a function in $L^\infty(G)$ given by $J(\varphi)(\cdot) = \la \lambda_\sigma(\cdot), \varphi \ra$.
By abuse of notation we will denote $J(\varphi)$ simply by $\varphi$.
This map allows us the following concrete description of $A(G,\sigma) = VN(G,\sigma)_*$.

\begin{prop}
We can identify $A(G,\sigma)$ as a subspace of $L^\infty(G)$ as follows.
\begin{align*}
A(G,\sigma)
& = \{ \varphi \in L^\infty(G) : \varphi(\cdot) =  \sum_{n\ge 1} \la \lambda_\sigma(\cdot)\xi_n, \eta_n \ra,\\
& \;\;\;\;\;\;\;\;\;\;  \;\;\;\;\;\;\;\;  \;\;\;\;\;\;\;\;  \;\; \xi_n, \eta_n \in L^2(G),\; \sum_{n\ge 1} \|\xi_n\| \cdot \|\eta_n\| <\infty\}
\end{align*}
with the norm
$$\|\varphi\|_{A(G,\sigma)} = \inf \left\{\sum_{n\ge 1} \|\xi_n\| \cdot \|\eta_n\| \right\},$$
where the infimum is taken over all possible such choice $\xi_n, \eta_n \in L^2(G)$, $\n\ge 1$.
With this identification the duality $(A(G,\sigma), VN(G,\sigma))$ becomes
$$\la \lambda_\sigma(\cdot), \varphi \ra = \varphi(\cdot),\;\; \varphi\in A(G,\sigma),$$
or equivalently
\begin{align}\label{eq-duality-twisted-Fourier}
\la \widetilde{\lambda_\sigma}(f), \varphi \ra = \int_G f(s) \varphi(s)ds,\;\; f\in L^1(G),\; \varphi\in A(G,\sigma).  
\end{align}
\end{prop}
	
\begin{cor}\label{cor-density-cpt-supp}
The space of functions in $A(G,\sigma)$ with compact support is dense in $A(G,\sigma)$.
\end{cor}

\begin{rem}
The fact that any element in $\varphi \in A(G)$ can be written as $\varphi(\cdot) = \la \lambda(\cdot)f, g\ra$ for some $f,g\in L^2(G)$ comes from the observation that $VN(G) \cong A(G)^*$ is in the standard form \cite{Tak2}. This still holds in the twisted case for the same reason that $VN(G,\sigma) \cong A(G,\sigma)^*$ is in the standard form, kindly communicated to us by Matt Daws. Thus, in the twisted case, we still have that, any element $\varphi \in A(G,\sigma)$ can we written as $\varphi(\cdot) = \la \lambda_\sigma(\cdot)f, g\ra$ for some $f,g\in L^2(G)$.
\end{rem}	

We can now see the relationship between $A(G,\sigma)$ and $B(G,\sigma)$.
\begin{prop}
The formal identity map
$$I:A(G,\sigma) \hookrightarrow B(G,\sigma)$$
is actually a complete isometric embedding.
\end{prop}
\begin{proof}
We only need to check the map $I$ is a complete isometry, equivalently the map $I^*: W^*(G,\sigma) \to VN(G, \sigma)$ is a complete quotient map. Note that $I^*$ is nothing but the sujective normal $*$-homomorphism guarranteed by the universality of $W^*(G,\sigma)$ and the left regular $\sigma$-representation $\lambda_\sigma: G\to VN(G,\sigma)$, which is known to be a complete quotient map by \cite[Lemma~III.2.2]{Tak1}.
\end{proof}

The Fourier algebra $A(G)$ is known be an ideal of the Fourier-Stieltjes algebra $B(G)$, which turns into the following with the twisting.
\begin{prop}\label{prop-bimodule}
For $\varphi \in B(G,\sigma)$ and $\psi \in A(G)$ we have $\varphi \psi \in A(G,\sigma)$.
The same conclusion holds for $\varphi \in A(G,\sigma)$ and $\psi \in B(G)$.
In particular, we have a complete contraction
$$m_\sigma: A(G,\sigma) \prt B(G) \to A(G,\sigma),\;\; f\otimes g \mapsto fg.$$
The above map makes $A(G,\sigma)$ an operator $B(G)$-bimodule and consequently an operator $A(G)$-bimodule. 	
\end{prop}
\begin{proof}
Let us focus on the first statement.
By polarization and a usual density argument it is enough to check the case $\psi \in C_c(G) \cap \mathcal{P}_1(G)$ and $\varphi \in \mathcal{P}_1(G, \sigma)$. 
Then, we know that $\varphi \psi \in \mathcal{P}_1(G, \sigma)$ with compact support, so that Proposition \ref{prop-cpt-supp-PD} gives us the conclusion that $\varphi \psi \in A(G,\sigma)$.
The second statement can be obtained similarly.
\end{proof}	
The adjoint map $\Gamma^\sigma = m^*_\sigma$ of a restriction of the above map $m_\sigma: A(G,\sigma) \prt A(G) \to A(G,\sigma),\;\; f\otimes g \mapsto fg$ is the {\bf twisted co-multiplication} given by
\begin{equation}\label{eq-twisted-co-product}
\Gamma^\sigma: VN(G,\sigma) \to VN(G,\sigma) \bar{\otimes} VN(G),\;\; \lambda_\sigma(s) \mapsto \lambda_\sigma(s) \otimes \lambda(s).
\end{equation}

We close this subsection with the description of the operator space structure on $A(G,\sigma) = VN(G,\sigma)_*$ by a norm formula for elements in $M_n(VN(G,\sigma))_*$. Recall two complete quotient maps
	$$\Phi: L^2(G)_R \prt L^2(G)_C \to VN(G,\sigma)_*,\; \xi\otimes \eta \mapsto \om_{\xi,\eta}$$
and
	$$\Psi: R_n \prt C_n \to (M_n)_*,\; e_i \otimes e_j \mapsto e_{ij}.$$
Then we get another complete quotient map
	$$\Phi \otimes \Psi: R_n \prt L^2(G)_R \prt L^2(G)_C \prt C_n \to M_n(VN(G,\sigma))_*\cong (M_n)_* \prt VN(G,\sigma)_*$$
after appropriate tensor flips. Via this map we have the following.
	\begin{prop} Every element of $M_n(VN(G,\sigma))_*$ is of the form
	\begin{equation}\label{eq-OS-twisted-Fourier}
	\varphi = \sum_{k\ge 1} \om_{(\xi^k_i)_i, (\eta^k_i)_i},
	\end{equation}
where $(\xi^k_i)^n_{i=1} \in R_n \prt L^2(G)_R$ and $(\eta^k_i)^n_{i=1} \in L^2(G)_C \prt C_n$ with
	$$\|\varphi\|_{M_n(VN(G,\sigma))_*} = \inf \Big\{\sum_{k\ge 1} \|(\xi^k_i)_i\|_{\ell^2_n(L^2(G))}\cdot \|(\eta^k_i)_i\|_{\ell^2_n(L^2(G))}\Big\},$$
	where the infimum runs over all possible expression in \eqref{eq-OS-twisted-Fourier}.
	\end{prop}

\section{Twisted multiplier spaces}\label{sec-twisted-multiplier}

In this section we introduce twisted multiplier spaces with respect to possibly different 2-cocycles on the group extending the discrete group case of B\'{e}dos/Conti \cite{BC2016-2}.

\begin{defn}
Let $\sigma_1, \sigma_2 \in Z^2(G,\tor)$ be locally continuous at the identity.
We say that a function $\varphi \in L^\infty(G)$ is a {\bf $(\sigma_1, \sigma_2)$-multiplier} if the operator
$$m_\varphi: A(G,\sigma_1) \to A(G,\sigma_2),\;\; \psi \mapsto \varphi \psi$$
is well-defined, i.e. $ \varphi \psi \in A(G,\sigma_2)$ for any $\psi \in A(G,\sigma_1)$.
We define two spaces of multipliers:
$$M(A(G,\sigma_1), A(G,\sigma_2)) := \{\varphi: G\to \Comp \;| \; \|m_\varphi\| <\infty\}$$
and
$$M_{cb}(A(G,\sigma_1), A(G,\sigma_2)) := \{\varphi: G\to \Comp \;| \; \|m_\varphi\|_{cb} <\infty\}.$$
When $\sigma_1 \equiv 1$ we simply write $M(A(G), A(G,\sigma_2))$ and $M_{cb}(A(G), A(G,\sigma_2))$. When $\sigma_1 = \sigma_2$ we simply write $MA(G,\sigma_1)$ and $M_{cb}A(G,\sigma_1)$.
\end{defn}
	
\begin{rem}\label{rem-multiplier-spaces-general}
Let $\sigma_1, \sigma_2, \sigma \in Z^2(G,\tor)$ be locally continuous at the identity.
\begin{enumerate}
\item By Proposition \ref{prop-bimodule} we have contractive inclusions
$$B(G,\sigma) \subseteq M_{cb}(A(G), A(G,\sigma)) \subseteq M(A(G), A(G,\sigma)).$$
		
\item We also have the contractive inclusion $M(A(G,\sigma_1), A(G,\sigma_2)) \subseteq L^\infty(G)$.
Indeed, taking adjoint of $m_\varphi$ we get another bounded map
$$M_\varphi = (m_\varphi)^* : VN(G,\sigma_2) \to VN(G,\sigma_1),\;\; \lambda_{\sigma_2}(s) \mapsto \varphi(s)\lambda_{\sigma_1}(s),$$
which clearly shows $\|\varphi\|_\infty \le \| M_\varphi  \|$.
\end{enumerate}		
\end{rem}

We have a characterization of	$M_{cb}(A(G,\sigma_1), A(G,\sigma_2))$ in the style of Gilbert \cite{Gil} and Jolissaint \cite{Joli92}.

\begin{thm}\label{thm-Jolissaint}
Let $\sigma_1, \sigma_2 \in Z^2(G,\tor)$ be locally continuous at the identity \color{black}  and $\varphi: G\to \Comp$ a function.
Then, $\varphi \in M_{cb}(A(G,\sigma_1), A(G,\sigma_2))$ if and only if there exists a  separable  Hilbert space $\Hi$ and bounded Borel measurable functions $\xi$ and $\eta$ from $G$ into $\Hi$ such that
\begin{equation}\label{eq-Gilbert-dec}
\sigma(s,t) \varphi(ts) = \la \xi(s), \eta(t) \ra,\;\; s,t\in G,
\end{equation}
where $\sigma = \overline{\sigma_1} \sigma_2\in Z^2(G,\tor)$. 	Moreover, we have
\begin{equation}\label{eq-cb-norm}
\|\varphi\|_{M_{cb}(A(G,\sigma_1), A(G,\sigma_2))} = \inf\{\|\xi\|_{L^\infty(G, \Hi)}\cdot \|\eta\|_{L^\infty(G, \Hi)}\},
	\end{equation}
	where the infimum runs over all possible choices of such $\xi(\cdot)$ and $\eta(\cdot)$.
	\end{thm}
\begin{proof}
Let $\varphi \in M_{cb}(A(G,\sigma_1), A(G,\sigma_2))$, which means the map $m_\varphi: A(G,\sigma_1) \to A(G,\sigma_2)$ is completely bounded. By taking the adjoint, we get another cb-map
	$$M_\varphi = (m_\varphi)^* : VN(G,\sigma_2) \to VN(G,\sigma_1),\;\; \lambda_{\sigma_2}(s) \mapsto \varphi(s)\lambda_{\sigma_1}(s).$$
By Wittstock's factorization theorem, we get a non-degenerate $*$-representation
	$$\Phi: C^*_r(G,\sigma_2) \to B(\Hi)$$
for some Hilbert space $\Hi$ and bounded maps $V_i: L^2(G) \to \Hi$, $i=1,2$ such that
	$$M_\varphi(\cdot) = V^*_2 \Phi(\cdot) V_1\;\; \text{and}\;\; \|M_\varphi\|_{cb} = \|V_1\| \cdot \|V_2\|.$$
Since $L^2(G)$ is separable, $\Hi$ can be replaced with its separable subspace 
$$\tilde{\Hi} := \overline{V_1(L^2(G))} \cup \overline{V_2(L^2(G))},$$
and $\Phi(\cdot)$ can then be replaced with $P_{\tilde{\Hi}}\Phi(\cdot)P_{\tilde{\Hi}}$, where $P_{\tilde{\Hi}}$ refers to the orthogonal projection from $\Hi$ onto $\tilde{\Hi}$. Thus, we may always assume $\Hi$ to be separable.
By composing the usual contractive embedding $L^1(G,\sigma_2) \hookrightarrow C^*_r(G,\sigma_2),\; f\mapsto \widetilde{\lambda_{\sigma_2}}(f)$ we get a $\sigma_2$-representation $\pi: G \to B(\Hi)$ such that $\tilde{\pi} = \Phi \circ \widetilde{\lambda_{\sigma_2}}$ thanks to Proposition \ref{prop-lifting-rep}. Thus, we have for $f\in L^1(G)$ that
	$$M_\varphi(\widetilde{\lambda_{\sigma_2}}(f)) = V^*_2 \tilde{\pi}(f) V_1.$$
Now we recall the bounded approximate identity	 $(\phi_U)_{U\in \mathcal{U}}$ of $L^1(G,\sigma_2)$ in Proposition \ref{prop-BAI}. Then, for $s\in G$ we have
	$$M_\varphi(\widetilde{\lambda_{\sigma_2}}(\lambda_{\sigma_2}(s)\phi_U)) = M_\varphi(\lambda_{\sigma_2}(s)\widetilde{\lambda_{\sigma_2}}(\phi_U)) \to M_\varphi(\lambda_{\sigma_2}(s))$$
in the weak$^*$-topology since $\widetilde{\lambda_{\sigma_2}}(\phi_U) \to 1_{VN(G,\sigma_2)}$ in the weak$^*$-topology as $U \to \{e\}$. On the other hand we have $\tilde{\pi}(\lambda_{\sigma_2}(s)\phi_U) = \pi(s)\tilde{\pi}(\phi_U)$. Note that for any $h,k \in \Hi$ we have
	$$\la \tilde{\pi}(\phi_U) h, k \ra = \int_G \psi_U(s) \la \pi(s)h, k \ra ds \to \la h, k\ra$$
as $U \to \{e\}$ since the function $\la \pi(\cdot)h, k \ra$ is bounded and locally continuous at the identity by Proposition \ref{prop-loc-continuity}. Thus, we have
$$V^*_2 \tilde{\pi}(\lambda_{\sigma_2}(s)\phi_U) V_1 \to V^*_2 \pi(s) V_1$$
in the weak operator topology on $B(\Hi)$. Consequently, we get
$$ \varphi(s)\lambda_{\sigma_1}(s) =  M_\varphi(\lambda_{\sigma_2}(s)) = V^*_2 \pi(s) V_1,\;\; s\in G.$$
	
Next, for a fixed unit vector $\xi_0\in L^2(G)$ we set
$$\xi(s) := \sigma(s,s^{-1})\pi(s^{-1})^*V_1\lambda_{\sigma_1}(s^{-1})\xi_0,\;\; \eta(t) := \pi(t)^*V_2\lambda_{\sigma_1}(t)\xi_0, \;\;s,t\in G.$$
Then, both $\xi(\cdot)$ and $\eta(\cdot)$ are bounded and  Borel measurable functions from $G$ to the separable Hilbert space $\Hi$.
Moreover, for any $s,t\in G$ we have
\begin{align*}
\la \xi(s), \eta(t) \ra
& = \sigma(s,s^{-1})\la \pi(s^{-1})^*V_1\lambda_{\sigma_1}(s^{-1})\xi_0, \pi(t)^*V_2\lambda_{\sigma_1}(t) \xi_0 \ra\\
& = \sigma(s,s^{-1})\la V^*_2 \pi(t)\pi(s^{-1})^*V_1\lambda_{\sigma_1}(s^{-1})\xi_0, \lambda_{\sigma_1}(t) \xi_0 \ra\\
& = \sigma(s,s^{-1}) \overline{\sigma_2(ts, s^{-1})}\la V^*_2 \pi(ts)V_1\lambda_{\sigma_1}(s^{-1})\xi_0, \lambda_{\sigma_1}(t) \xi_0 \ra\\
& = \sigma(s,s^{-1}) \overline{\sigma_2(ts, s^{-1})}\varphi(ts)\la \lambda_{\sigma_1}(ts)\lambda_{\sigma_1}(s^{-1})\xi_0, \lambda_{\sigma_1}(t) \xi_0 \ra\\
& = \sigma(s,s^{-1}) \overline{\sigma_2(ts, s^{-1})}\varphi(ts)\sigma_1(ts,s^{-1})\la \lambda_{\sigma_1}(t)\xi_0, \lambda_{\sigma_1}(t) \xi_0 \ra\\
& = \sigma(t,s)\varphi(ts). 
\end{align*}
Finally, we have
	\begin{equation}\label{eq-cb-norm-upper}
	\|\xi\|_{L^\infty(G, \Hi)}\cdot \|\eta\|_{L^\infty(G, \Hi)} \le \|V_1\| \cdot \|V_2\| = \|M_\varphi\|_{cb}.
	\end{equation}
	
Conversely, we assume that \eqref{eq-Gilbert-dec} holds.	Consider a function $\psi \in A(G,\sigma)$ of the form $\psi(\cdot) = \la \lambda_{\sigma_1}(s) f, g \ra$ for $f,g\in L^2(G)$. Then, we have
\begin{align*}
\varphi(s)\psi(s)
& = \int_G \varphi(s) \sigma_1(s,s^{-1}t) f(s^{-1}t)\overline{g(t)}dt\\
& = \int_G \sigma_2(s,s^{-1}t) \overline{\sigma(s^{-1}t, t^{-1}s)} \la \xi(t^{-1}s), \eta(t) \ra f(s^{-1}t)\overline{g(t)}dt\\
& = \int_G \sigma_2(s,s^{-1}t) \overline{\sigma(s^{-1}t, t^{-1}s)} \sum_{i\in \mathbb{N}}\la \xi(t^{-1}s), e_i \ra \cdot  \la e_i, \eta(t) \ra f(s^{-1}t)\overline{g(t)}dt\\
& = \sum_{i\in \mathbb{N}} \int_G \sigma_2(s,s^{-1}t) \overline{\sigma(s^{-1}t, t^{-1}s)}\la \xi(t^{-1}s), e_i \ra f(s^{-1}t)\overline{\la \eta(t), e_i\ra g(t)}dt\\
& = \sum_{i\in \mathbb{N}}\la \lambda_{\sigma_2}(s) h_i, k_i \ra,
\end{align*}
where $(e_i)_{i\in \mathbb{N}}$ is a fixed orthonormal basis of $\Hi$ and $h_i, k_i \in L^2(G)$ are given by
$$h_i(s) :=\overline{\sigma(s,s^{-1})\la \xi(s^{-1}), e_i \ra f(s)},\;\; k_i(t) := \la \eta(t), e_i \ra g(t),\; s,t\in G.$$
Then, we have $\varphi\psi \in A(G,\sigma)$ with
$$\|\varphi\psi\|_{A(G,\sigma)}\le \sum_{i\in \mathbb{N}}\|h_i\| \cdot \|k_i\| \le (\sum_{i\in \mathbb{N}} \|h_i\|^2)^{\frac{1}{2}} (\sum_{i\in \mathbb{N}}\|k_i\|^2)^{\frac{1}{2}}$$ and
	\begin{align*}
		\|\varphi\psi\|_{A(G,\sigma)}
		& \le \sum_{i\in \mathbb{N}}\|h_i\| \cdot \|k_i\|
		\le (\sum_{i\in \mathbb{N}}\|h_i\|^2)^{\frac{1}{2}} (\sum_{i\in \mathbb{N}}\|k_i\|^2)^{\frac{1}{2}}\\
		& = (\sum_{i\in \mathbb{N}}\int_G |\la \xi(s^{-1}), e_i \ra f(s)|^2ds)^{\frac{1}{2}} (\sum_{i\in \mathbb{N}}\int_G |\la \eta(t), e_i \ra g(t)|^2dt)^{\frac{1}{2}}\\
		& \le \|\xi\|_{L^\infty(G, \Hi)}\cdot \|\eta\|_{L^\infty(G, \Hi)}\cdot \|f\|_2 \cdot \|g\|_2.
	\end{align*}

 \color{black}

A general element $\psi \in A(G,\sigma)$ is actually of the form $\psi(\cdot) = \sum_{n\ge 1}\la \lambda_{\sigma_1}(s) f_n, g_n \ra$ for $f_n, g_n\in L^2(G)$. Then the above argument tells us that we have
	$$\|\varphi\psi\|_{A(G,\sigma)} \le  \|\xi\|_{L^\infty(G, \Hi)}\cdot \|\eta\|_{L^\infty(G, \Hi)}\cdot \sum_{n\ge 1}\|f_n\|_2 \cdot \|g_n\|_2,$$
so that we have
	$$\|\varphi\psi\|_{A(G,\sigma)} \le  \|\xi\|_{L^\infty(G, \Hi)}\cdot \|\eta\|_{L^\infty(G, \Hi)}\cdot \|\psi\|_{A(G,\sigma)}$$
by taking infimum over all such choices of $f_n, g_n \in L^2(G)$. This explains that the multiplier $m_\varphi$ is bounded. The exactly same argument still works in the matrix valued case using the description of \eqref{eq-OS-twisted-Fourier}, so that we actually have that $m_\varphi$ is completely bounded with
	$$\|m_\varphi\|_{cb} \le \|\xi\|_{L^\infty(G, \Hi)}\cdot \|\eta\|_{L^\infty(G, \Hi)}.$$
Together with \eqref{eq-cb-norm-upper} we get \eqref{eq-cb-norm}.
\end{proof}

\begin{cor}
Let $\sigma_1, \sigma_2 \in Z^2(G,\tor)$ be locally continuous at the identity. 
The multiplier space $M_{cb}(A(G,\sigma_1), A(G,\sigma_2))$ depends only on the difference of 2-cocycles, namely $\overline{\sigma_1} \sigma_2$, so that we have
$$M_{cb}(A(G,\sigma_1), A(G,\sigma_2)) =  M_{cb}(A(G), A(G,\overline{\sigma_1} \sigma_2))$$
with the same norms.
In particular, for any $\sigma \in Z^2(G,\tor)$ locally continuous at the identity, we have
$$M_{cb}A(G) = M_{cb}A(G,\sigma)$$
with the same norms.	
\end{cor}

	\begin{rem}
		\begin{enumerate}
			\item The conclusion for $\sigma_1 = \sigma_2$ first appeared in \cite{BC2009} for discrete groups.
			
			\item For the special case $G=\mathbb{Z}^{2d}$ with the 2-cocycle in Example \ref{ex-Q-torus-plane} the above corollary and its $L^p$-counterpart ($1\le p \le \infty$) were studied in \cite{CXY2012}. However, it has been shown that there are $L^p$-bounded ($1\le p \le \infty$) Fourier multipliers on the non-commutative torus, which are not completely bounded on $L^p$. See \cite{Ricard2016} for the details.
			
			\item For $G=\mathbb{R}^{2d}$ the above corollary says that the completely bounded $L^p$-Fourier multipliers on quantum Euclidean spaces are the same as those on classical Euclidean spaces for $p = 1$ or $\infty$.
		\end{enumerate}
	\end{rem}

The following proposition follows immediately from Theorem \ref{thm-Jolissaint} and the definition of $B(G,\sigma)$.
\begin{prop}\label{prop:inclusion bg mag}
Let $\sigma, \sigma_1, \sigma_2 \in Z^2(G,\tor)$ be locally continuous at the identity, which satisfy $\sigma=\overline{\sigma_1 } \sigma_2$.
Then, we have
$$B(G,\sigma) \subseteq  M_{cb}\big(A(G,\sigma_1),A(G,\sigma_2)\big)\subseteq M\big(A(G,\sigma_1),A(G,\sigma_2)\big).$$
\end{prop}
	
\begin{rem}
The idea of using two different $2$-cocycles for the multiplier spaces is crucial in our investigation, which is partially due to the fact that $A(G,\sigma)$ is merely a $A(G)$-bimodule with respect to pointwise multiplication. In \cite{Ricard2016} a similar idea has been already demonstrated in the special case of non-commutative torus $C(\tor^2_\theta)$, namely using two different irrational parameter $\theta$ for the study of associated Fourier multipliers.
\end{rem}

	
\section{Amenability and the twisted multiplier space}\label{sec-amen}


Let us fix $\sigma \in Z^2(G,\tor)$, which is locally continuous at the identity.
We will see in this section that the equality $B(G,\sigma) = M_{cb}(A(G),A(G,\sigma))$, and consequently the equality $B(G,\sigma) = M(A(G),A(G,\sigma))$ characterizes the amenability of the group as in the untwisted case.
The following is one ingredient we need here, which is the easy direction of a twisted version of Hulanicki's theorem \cite{Hul64}.

\begin{prop}\label{amenable-C-Cr}
If $G$ is amenable, then the canonical quotient map
	$$Q_\sigma :C^*(G,\sigma)\rightarrow C_r^*(G,\sigma),\;\; \widetilde{\pi^\sigma_u}(f) \mapsto \widetilde{\lambda_\sigma}(f),\;\; f\in L^1(G)$$
is an isometric isomorphism.
\end{prop}

\begin{proof}
From the density argument it is enough to check the following inequality for any $f\in C_c(G)$:
\begin{equation}\label{pi-leq-left}
\|\widetilde{\pi^\sigma_u}(f)\|_{C^*(G,\sigma)} \le \|\widetilde{\lambda_\sigma}(f)\| _{C^*_r(G,\sigma)} .
\end{equation}
Since 
$$\|\widetilde{\pi^\sigma_u}(f)\|_{C^*(G,\sigma)}  ^2 = \|\widetilde{\pi^\sigma_u}(f)^*\widetilde{\pi^\sigma_u}(f)\| _{C^*(G,\sigma)}  = \|\widetilde{\pi^\sigma_u}(f^\star *_\sigma f)\| _{C^*(G,\sigma)}  ,$$
 by Proposition \ref{prop-P1-equiv}, we can find $\varphi \in \mathcal{P}_1(G,\sigma)$ such that 
 $$\|\widetilde{\pi^\sigma_u}(f)\|_{C^*(G,\sigma)}  ^2 = \int_G f^\star *_\sigma f(s) \varphi(s) ds  .$$
Since $G$ is amenable, there exists a net $(\psi_\alpha)\subset \mathcal{P}_1(G)$ with compact support such that $\psi_\alpha \to 1$ pointwise, so that $\psi_\alpha(s) \varphi(s) \to \varphi(s)$ uniformly on the support of $f^\star *_\sigma f$.
Each function $\psi_\alpha \varphi$ has compact support, and by Corollary \ref{cor-sigma-pd}, belongs to $\mathcal{P}_1(G,\sigma)$, so that we appeal to Proposition \ref{prop-cpt-supp-PD} to write $\psi_\alpha(\cdot) \varphi(\cdot)  = \langle \lambda_\sigma(\cdot) \eta_\alpha, \eta_\alpha\rangle$
for some unit vector $\eta_\alpha\in L^2(G)$.
Now we have
\begin{align*}
\|\widetilde{\pi^\sigma_u}(f)\|_{C^*(G,\sigma)}  ^2 
& = \lim_{\alpha} \int_G f^\sigma *_\sigma f (s) \langle \lambda_\sigma(s) \eta_\alpha, \eta_\alpha\rangle\, ds\\
&= \lim_{\alpha} \langle \widetilde{ \lambda_\sigma }(f^\sigma *_\sigma f ) \eta_\alpha, \eta_\alpha \rangle\\
& \leq \limsup_\alpha\| \widetilde{ \lambda_\sigma}(f) \|^2_{C_r^*(G,\sigma)} \|\eta_\alpha\|_2^2\\
& = \| \widetilde{ \lambda_\sigma}(f) \|^2_{C_r^*(G,\sigma)} ,
\end{align*}
completing the proof.
\end{proof}

	\begin{rem}

	The converse direction of the above proposition, which is the hard part of a twisted version of Hulanicki's theorem \cite{Hul64} is not yet known.
	 \end{rem}

Now we go back to the study of multipliers, and consider the relation between the amenability of $G$ and the following inclusions from Proposition \ref{prop:inclusion bg mag}
	$$B(G,\sigma) \subseteq M_{cb}\big(A(G),A(G,\sigma)\big) \subseteq M\big(A(G),A(G,\sigma)\big).$$
First, with aide of the above proposition, we deduce the following.

\begin{thm}\label{amen-imply-eq}
If the locally compact group $G$ is amenable, then we have
	$$B(G,\sigma) = M_{cb}\big(A(G),A(G,\sigma)\big) = M\big(A(G),A(G,\sigma)\big).$$	
\end{thm}
\begin{proof}
For $\psi \in M\big(A(G),A(G,\sigma)\big) \subseteq  L^\infty(G)$ (by Remark \ref{rem-multiplier-spaces-general}) we consider the following composition.
\[\begin{array}{cccccccc}
		\Psi: &C^*_r(G,\sigma) &\stackrel{M_\psi}{\longrightarrow}& C^*_r(G) &\cong & C^*(G) &\stackrel{\varphi}{\longrightarrow} &  \mathbb{C}\\
		&\widetilde{\lambda_\sigma}(f) = \int_G f(s)\lambda_\sigma(s)ds & \mapsto & \widetilde{\lambda}(f\psi) &\mapsto&\widetilde{\pi_u}(f\psi) & \mapsto & \int_G f(s)\psi(s)ds
		\end{array}\]
		Here, we have $C^*_r(G) \cong C^*(G)$ by amenability and $\varphi$ denotes the trivial homomorphism from $C^*(G)$ into $\Comp$. We can see that the functional $\Psi  \in C^*_r(G,\sigma)^* \cong C^*(G,\sigma)^*$ by Proposition \ref{amenable-C-Cr}, which means that $\psi \in B(G,\sigma)$.
\end{proof}	

\begin{rem}
If the locally compact group $G$ is amenable and the 2-cocycle $\sigma \in Z^2(G,\tor)$ is locally continuous at the identity, we have already seen that $B(G)=  M_{cb}A(G) =  M_{cb}A(G, \sigma)$. However, the inclusion $M_{cb}A(G, \sigma)\subseteq MA(G, \sigma)$ could be proper even for abelian groups, which also implies that the measure algebra $ M(\widehat G)=M_{cb}A(G ) = MA(G) $ is properly contained in $ MA(G, \sigma)$.
See \cite[Remark~7.7]{CXY2012} for $G= \mathbb{Z}^\infty$, and \cite[Proposition~4.1]{Ricard2016} for $G= \mathbb{Z}^2$.
\end{rem}

Now we consider the converse of Theorem \ref{amen-imply-eq}, which is the hard part of the characterization of amenability through twisted multiplier spaces.

\begin{thm}\label{eq-imply-amen}
If $B(G,\sigma) = M_{cb}\big(A(G),A(G,\sigma)\big)$, then $G$ is amenable.
\end{thm}
The proof of Theorem \ref{eq-imply-amen} is much subtler; we will divide its proof into two cases, given in the next two subsections seperately.

\subsection{Proof of Theorem \ref{eq-imply-amen}: discrete case}\label{subsec-discrete}

In this section, we will always assume that $G$ is a discrete group, and prove Theorem \ref{eq-imply-amen} in this case. We will adapt the argument of Bo\.zejko \cite{Boz85} to the twisted case; a crucial observation is that the Littlewood function space $T_2$ below (introduced in \cite{Var1974}) comes from another function space $t_2$, which is invariant under the multiplication with respect to a function with values in $\tor$ (e.g. $2$-cocycle). Before the proof we need some preparations.

Let us recall a classical result of Grothendieck on Schur multipliers.

\begin{thm} (\cite[Theorem 5.1]{Pisier1991}) A function $f: G\times G \to \Comp$ is a (completely) bounded Schur multiplier on $B(\ell^2(G))$ if and only if there exist a Hilbert space $\Hi$ and functions $\xi$ and $\eta$ from $G$ into $\Hi$ such that
		\begin{equation}\label{eq-Grothendieck-dec}
		f(t,s) = \la \xi(s), \eta(t) \ra,\;\; s,t\in G
		\end{equation}
	Moreover, we have
		$$	\|M_f\|_{cb} = \|M_f\| = \inf\{\|\xi\|_{\ell^\infty(G, \Hi)}\cdot \|\eta\|_{\ell^\infty(G, \Hi)}\},$$
	where the infimum runs over all possible choices of such $\xi(\cdot)$ and $\eta(\cdot)$. Here, $M_f$ is the map
		$$M_f: B(\ell^2(G)) \to B(\ell^2(G)),\;\; A=(A(s,t))_{s,t\in G}\mapsto (f(s,t)A(s,t))_{s,t\in G}.$$
\end{thm}

\begin{rem}\label{rem-recast-Jolissaint}
Theorem \ref{thm-Jolissaint} can be restated as follows: a function $\varphi:G \to \Comp$ is a completely bounded $(1,\sigma)$-multiplier if and only if $f(s,t) = \sigma(s,t)\varphi(ts)$ is a (completely) bounded Schur multiplier on $B(\ell^2(G))$ with $\|\varphi\|_{M_{cb}(A(G),A(G,\sigma))} = \|M_f\| = \|M_f\|_{cb}$. 
\end{rem}

Now we recall some function spaces.

\begin{defn}\label{def-T2}
We define
	$$X_1 :=\big\{  \psi: G\times G\rightarrow \mathbb{C}: \|\psi\|_{X_1}=\sup_{s\in G} (\sum_{t\in G}|\psi(s,t)|^2)^{\frac{1}{2}}< \infty   \big\},$$
	$$X_2 :=\big\{ \psi: G\times G\rightarrow \mathbb{C}: \|\psi\|_{X_2}=\sup_{t\in G}(\sum_{s\in G}|\psi(s,t)|^2)^{\frac{1}{2}}< \infty   \big\}$$
and
	$$t_2 := X_1+X_2 = \big\{ \psi: G\times G\rightarrow \mathbb{C} : \exists \,\psi_1\in X_1, \psi_2\in X_2\; \text{such that}\; \psi=\psi_1+\psi_2
\big\}$$
	with the norm $\displaystyle \|\psi\|_{t_2}=\inf_{\psi=\psi_1+\psi_2} \left\{\|\psi_1\|_{X_1} + \|\psi_2\|_{X_2}\right\}.$
	
	We also define the space $T_2$ of Littlewood functions on $G$ as follows.
	\[T_2 = T_2(G) :=\{\varphi:G\rightarrow \mathbb{C}: \; f(s,t):=\varphi(st)\in t_2, s,t\in G\}\] 
		with the norm $\|\varphi\|_{T_2}:= \|f\|_{t_2}$.
\end{defn}

\begin{rem}\label{rem-T2}
\begin{enumerate}
	\item We may define $T_p$ for $1\le p <\infty$ similarly, but we focus only on $T_2$ in this article.
	
	\item Letting $\psi_1(s,t) = \varphi(st)$, $\psi_2\equiv 0$ we can easily see that $\ell^2(G) \subseteq T_2$, contractively.
	
	\item From the definition of $t_2$ it is straightforward to check that for any $\psi\in t_2$ and $f:G\times G \to \tor$ we have $\psi f \in t_2$ with $\|\psi f \|_{t_2}=\|\psi \|_{t_2}$. Indeed, for any decomposition $\psi=\psi_1+\psi_2$ we get the corresponding decomposition $\psi f = \psi_1 f + \psi_2 f$, so that $\|\psi f \|_{t_2}\le \|\psi_1 f\|_{X_1} + \|\psi_2 f\|_{X_2} = \|\psi_1 \|_{X_1} + \|\psi_2 \|_{X_2}$. Taking infimum over all such decompositions we get $\|\psi f \|_{t_2}\le \|\psi \|_{t_2}$. Now we get the desired conclusion by repeating the same argument for $\psi f$ and $\bar{f}$.
	
	\item It is also straightforward to check that any $\psi \in t_2$ is a (completely) bounded Schur multiplier with
		$$\|M_\psi\| \le \|\psi\|_{t_2}.$$
	Indeed, for any decomposition $\psi=\psi_1+\psi_2$ with $\psi_1\in X_1$ and $\psi_2\in X_2$ we set $\xi(s) = \sum_{s'\in G}\psi_1(s,s')\delta_{s'} \in \ell^2(G)$ and $\eta(t) = \delta_t \in \ell^2(G)$, $s,t\in G$, where $\delta_t$ is the delta function on $t\in G$.
Then we have $\psi_1(s,t) = \la \xi(s), \eta(t) \ra$, $s,t\in G$ and so that $\|M_{\psi_1}\| \le \sup_{s\in G}\|\xi(s)\|_2 \cdot \sup_{t\in G}\|\eta(s)\|_2 = \|\psi_1\|_{X_1}.$
Similarly we have $\|M_{\psi_2}\| \le \|\psi_2\|_{X_2}$, which means that $\|M_\psi\| \le \|\psi_1\|_{X_1} + \|\psi_2\|_{X_2}$. Taking infimum over all such decompositions we get $\|M_\psi\| \le \|\psi\|_{t_2}.$	
\end{enumerate}
\end{rem}

\begin{lem}
Let $M(\sigma) := \{\varphi: G\to \Comp\, |\, \|\varphi \|_{M(\sigma)} <\infty\}$ the space of all bounded multipliers from $\ell^ \infty(G)$ to $M_{cb}(A(G), A(G,\sigma))$, where
	$$\|\varphi \|_{M(\sigma)} := \|M_\varphi: \ell^\infty(G) \to M_{cb}(A(G), A(G,\sigma)),\; g\mapsto \varphi g \|.$$
Then we have the following contractive inclusion.
	\[T_2\subseteq M(\sigma).\]
\end{lem}
\begin{proof}
Fix $\varphi \in T_2$. From a standard extreme point argument it is enough to check that $1_E \varphi \in M_{cb}(A(G), A(G,\sigma))$ with $\|1_E \varphi\|_{M_{cb}(A(G), A(G,\sigma))}\le \|\varphi\|_{T_2}$ for any $E\subseteq G$. Indeed, for the decomposition $1_E=\frac{1}{2}[(1_E-1_{G\setminus E})+ 1_G] = \frac{1}{2}(\xi_E+ 1_G)$ we can see that the function $(s,t)\mapsto \varphi(st)$ is in $t_2$, so that the function $(s,t)\mapsto \sigma(t,s)\xi_E(st)\varphi(st)$ is a bounded Schur multiplier on $B(\ell^2(G))$. By Theorem \ref{thm-Jolissaint}, Remark \ref{rem-recast-Jolissaint} and Remark \ref{rem-T2} (4) we know that $\xi_E \varphi \in M_{cb}(A(G), A(G,\sigma))$ with $\|\xi_E \varphi \|_{M_{cb}(A(G), A(G,\sigma))}\le \|\varphi\|_{T_2}$. Replacing $\xi_E$ by $1_G$, we have also $\varphi \in M_{cb}(A(G), A(G,\sigma))$  with
	$$\|\varphi \|_{M_{cb}(A(G), A(G,\sigma))}\le \|\varphi\|_{T_2},$$
which, in turn, tells us that $1_E \varphi \in M_{cb}(A(G), A(G,\sigma))$ with
	$$\|1_E \varphi \|_{M_{cb}(A(G), A(G,\sigma))}\le \|\varphi\|_{T_2}.$$
		\end{proof}

Now we are going to complete the proof of Theorem \ref{eq-imply-amen} when the group $G$ is discrete.

\begin{proof}[Proof of Theorem \ref{eq-imply-amen}]

Recall that a Banach space $X$ is said to be of {\it cotype 2}, if there exists a universal constant $C$, such that 
	\[\left(\sum_{i=1}^n\|x_i\|^2\right)^{\frac{1}{2}}\leq C\left(\int_0^1\|\sum_{i=1}^n r_i(t)x_i\|^2 dt\right)^{\frac{1}{2}}\]
	for all finite $\{x_i\}_{i=1}^n\subseteq X$, where $\{r_i\}_{i=1}^\infty$ is the sequence of Rademacher functions on $[0,1]$. From the assumption we have $M_{cb}(A(G), A(G,\sigma))=B(G,\sigma)$ the dual of a $C^*$-algebra $C^*(G,\sigma)$, so that it is of cotype 2 as a Banach space (\cite{TJ1974, Pisier1978}). For any finite subset $\{s_1, \cdots, s_n\} \subseteq G$ and $(a_j)^n_{j=1}\subseteq \Comp$ we consider $g=\sum_{n=1}^n a_j \delta_{s_j} : G \to \Comp$. For $\varepsilon \in [0,1]$ we set $g_\varepsilon := \sum_{n=1}^n a_j r_j(\varepsilon) \delta_{s_j}: G \to \Comp$.
From the cotype $2$ condition we have a universal constant $C$ such that,
	$$C (\sum^n_{j=1}|a_j|^2)^{\frac{1}{2}} \le \int_0^1 \|g_\varepsilon(s)\|_{M_{cb}(A(G), A(G,\sigma))}d\varepsilon \le \|g\|_{M(\sigma)}.$$
Since the choice of $\{s_1, \cdots, s_n\}$ and $(a_j)^n_{j=1}$ was arbitrary we have $M(\sigma)\subseteq \ell^2(G)$, which in turn gives us $T_2 \subseteq M(\sigma) \subseteq \ell^2(G) \subseteq T_2$ and consequently $T_2 = \ell^2(G)$. This implies that $G$ is amenable by a result of J. Wysoczanski \cite{Wys88} (or \cite[Theorem 2.5]{Pisier1991}).
\end{proof}

\subsection{Proof of Theorem \ref{eq-imply-amen}: non-discrete case}\label{subsec-non-discrete}

This subsection deals with the proof of Theorem \ref{eq-imply-amen} when $G$ is a non-discrete locally compact group, which is the most technical part in this paper. We keep the assumption that $\sigma \in Z^2(G,\tor)$ is locally continuous at the identity.
Our treatment follows Losert \cite{Losert} rather closely with modifications to the twisted case and to the operator space level.

We first introduce the sufficient condition used in \cite{Losert}, for the amenability of the group $G$.
It is based on a conclusion on the space of left-uniformly continuous functions $LUC(G)$. If $M$ is a mean on $L^\infty (G) $ or $LUC(G)$ and $s\in G$, we write
	$$d(M,x) =\sup \{  |M(\lambda(x) f -f)| : \|f\|_\infty \le 1  \}.$$

\begin{prop}[\cite{Losert}]\label{prop-LUC-mean}
Assume there exists $c<2$ such that for each $x_1,\cdots, x_N \in G$ (not necessarily distinct points) there exists a mean $M$ on $LUC(G)$ with $\sum_{k=1}^N d(M,x_k) \leq c N $. Then $G$ is amenable.
\end{prop}

Before we proceed to the details of the proof of Theorem \ref{eq-imply-amen}, we would like to introduce two different types of products between $\lambda_\sigma(s)$ and $u\in A(G,\sigma)$ as follows.
\begin{defn}
For $u\in A(G,\sigma)$ and $T\in VN(G,\sigma)$ we define $T\cdot u\in VN(G)$ by
$$\la T \cdot u, v \ra := \la T, uv \ra,\;\; v\in A(G).$$
For $s\in G$ the above defines an operator $\lambda_\sigma(s)\cdot u\in VN(G)$.
We also define
$$\lambda_\sigma(s) \bullet u := \overline{\lambda_\sigma(s)(\bar{u})} \in L^\infty(G),$$
where $\lambda_\sigma(s)$ is regarded as an operator acting on $L^\infty(G)$.
\end{defn}

\begin{rem}
\begin{enumerate}
\item For $u\in A(G,\sigma)$ we can easily see that $T\cdot u = m^*_u(T)$, where $m_u: A(G)\to A(G,\sigma)$ is the multiplier associated to $u$.
Moreover, from \eqref{eq-duality-twisted-Fourier}, we can easily see the following.
$$\widetilde{\lambda_\sigma}(f)\cdot u = \widetilde{\lambda_\sigma}(uf) \in C^*_r(G,\sigma),\;\; f\in L^1(G).$$

\item When $\sigma$ is continuous, we have an alternative pointwise description, namely for $s\in G$
$$\lambda_\sigma(s) \cdot u = u(s)\lambda_\sigma(s).$$
\end{enumerate}
\end{rem}

\begin{prop}\label{prop-two-actions}
For $s\in G$, $u\in A(G,\sigma)$ and $T\in VN(G,\sigma)$ we have the following.
\begin{enumerate}
\item $\lambda_\sigma(s) \bullet u \in A(G,\sigma)$ and
\begin{equation}\label{eq-bullet-dual}
\la T, \lambda_\sigma(s) \bullet u \ra = \la \lambda_\sigma(s)^*T, u \ra.
\end{equation}
			
\item The map 
$$\Phi_s : A(G,\sigma) \to A(G,\sigma),\; f\mapsto \lambda_\sigma(s)\bullet f$$
is a completely isometric isomorphism.	
\end{enumerate}
\end{prop}

\begin{proof}
(1) We begin with $u\in A(G,\sigma)$ of the form $u(\cdot) =\langle \lambda_\sigma(\cdot) \xi, \eta\rangle$ for some $\xi,\eta \in L^2(G)$.
Then, we have 
\begin{equation}\label{bullet-inner}
\lambda_\sigma(s) \bullet u(\cdot) =\langle \lambda_\sigma(\cdot) \xi, \lambda_\sigma(s)\eta \rangle
\end{equation}
by direct computations.
This explains $\lambda_\sigma(s) \bullet u\in A(G,\sigma)$ with $\|\lambda_\sigma(s) \bullet u\|_{A(G,\sigma)} \le \|\xi\|\cdot\|\eta\|$,
and the general case can be obtained similarly.

The second statement is also straightforward to the case $T = \lambda_\sigma(t)$, $t\in G$ and we again appeal to linearity and density for the general case.

\vspace{0.5cm}

(2) By \eqref{bullet-inner}, we know that $\Phi_s$ is a contraction. On the other hand, we can readily see that 
$$\lambda_\sigma(s) ^* \bullet (\lambda_\sigma(s) \bullet u)  = \overline{\sigma(s,s^{-1})}^2 u, $$
or equivalently, $u =\sigma(s,s^{-1})  \lambda_\sigma(s^{-1}) \bullet(\lambda_\sigma(s) \bullet u)$. So $\|u\|_{A(G,\sigma) }\le   \| \lambda_\sigma(s) \bullet u \|_{A(G,\sigma)}$ by the contractivity of $\Phi_s$. Thus $\Phi_s$ is an onto isometry.

For the operator space level, we recall the characterization of canonical elements in $[\fM_n(VN(G,\sigma))]_*$ given in \eqref{eq-OS-twisted-Fourier}, from which we deduce 
$$I_n\otimes \Phi_s \Big( \sum_k \big\langle  I_n \otimes \lambda_\sigma(\cdot) (\xi_i^k)_i, (\eta _j ^k) _j    \big\rangle   \Big) =\sum_k \big\langle  I_n \otimes \lambda_\sigma(\cdot)(\xi_i^k)_i, I_n\otimes \lambda_\sigma(s) (\eta _j ^k) _j    \big\rangle.$$
Hence $\Phi_s $ is a complete contraction. Using the same trick as in the Banach space level, we can finally see that $\Phi_s$ is an onto complete isometry for every $s\in G$.
\end{proof}

With the help of above actions, we can get the following technical lemma and induction proposition for the proof of Theorem \ref{eq-imply-amen}.
In this lemma and proposition, $f$ denotes a fixed continuous function on $\Real$ such that
$0 \leq f(t) \leq 1$, $f(t) = 0$ for $t \leq \frac 1  4$ and $f(t) = 1$ for $t\geq \frac 1 2 $.

\begin{lem}\label{lem-inner-control}
Let $T\in VN(G,\sigma)$ with $0\leq T\leq 1$ and $u\in \mathcal{P}_1(G,\sigma)$ with compact support. Suppose that $\langle T, u \rangle > 1-\frac{\delta^2}{2}$.
\begin{enumerate}
\item If $\langle \lambda_\sigma(s)^* T \lambda_\sigma(s) , u \rangle > 1-\frac{\delta^2}{2}$, then we have 
$$\langle f(T) \lambda_\sigma(s) f(T) ,  \lambda_\sigma(s)\bullet u \rangle >1-2\delta.$$
\item Suppose that $(T_i)^n_{i=1}$ is a sequence of commuting operators in $VN(G,\sigma)$ such that $0\le T_i \le 1$ and $\langle \sum^n_{i=1}T_i^2, u \rangle <\big(\frac{\delta^2}{4n}\big)^2$ for $1\leq i\leq n$, then 
$$\la \Big(\prod^n_{i=1}(1-T_i )\Big) T \Big(\prod^n_{i=1}(1-T_i )\Big), u \rangle > 1-\delta^2.$$
\end{enumerate}
\end{lem}
\begin{proof}
Note that Proposition \ref{prop-cpt-supp-PD} ensures that $u(\cdot) = \la \lambda_\sigma(\cdot)h, h \ra$ for some unit vector $h\in L^2(G)$. Then, we may repeat the the same argument of \cite[Lemma 2]{Losert} for the rest of the proof with the help of \eqref{eq-bullet-dual}.
\end{proof}

The following is a twisted version of \cite[Lemma 1]{Losert}, where the same argument works.
\begin{lem}\label{lem-cpt}
For $T \in C^*_r(G,\sigma)$ and a compact subset $K\subseteq G$, the restricted operator $T|_{L^2(K)}$ is compact.
\end{lem}

\begin{prop}\label{prop-ell_1}
Assume that $G$ is non-discrete. For any $(u_i)^n_{i=1} \subseteq  \mathcal{P}_1(G,\sigma)$ with compact support, 
$(s_i)^n_{i=1} \subseteq G$ and $\eps>0$ there are $(R_i)^n_{i=1}, (\overline{R}_i)^n_{i=1} \subseteq C^*_r(G,\sigma)$ and $(w_i)^n_{i=1} \subseteq \mathcal{P}_1(G)$ such that
\begin{enumerate}
\item $0\le R_i, \overline{R}_i \le 1$, $1\le i\le n$;

\item $\overline{R}_i\overline{R}_j = \overline{R}_j\overline{R}_i = 0$, $1\le i\ne j\le n$;

\item $R_i\overline{R}_i = \overline{R}_iR_i = R_i$, $1\le i \le n$;

\item $\la R_i\lambda_\sigma(s_i)R_i, \lambda_\sigma(s_i) \bullet (u_iw_i) \ra > 1 - \eps$.
\end{enumerate}
In particular, for any $(\mu_i)^n_{i=1} \subseteq \Comp$ we have
$$\|\sum^n_{i=1}\mu_i \lambda_\sigma(s_i) \bullet (u_iw_i)\|_{A(G,\sigma)} > (1-\eps)\sum^n_{i=1}|\mu_i|.$$
\end{prop}
\begin{proof}
We begin with the case $n=1$. Proposition \ref{prop-cpt-supp-PD} says that $u_1(\cdot) = \la \lambda_\sigma(\cdot)h, h \ra$ for some unit vector $h\in L^2(G)$. Thus, we have $\la T, u_1 \ra = \la Th, h \ra$ for any $T \in C^*_r(G,\sigma)$. By Kaplansky density theorem there is $R\in C^*_r(G,\sigma)$ such that $0\le R \le 1$ and
	$$\la R, u_1 \ra,\; \la\lambda_\sigma(s_1)^* R \lambda_\sigma(s_1), u_1 \ra > 1 - \frac{\eps^2}{8}.$$
Then, Lemma \ref{lem-inner-control} (1) tells us that $\la f(R) \lambda_\sigma(s_1) f(R), \lambda_\sigma(s_1) \bullet u_1 \ra > 1 - \eps$. Thus, we get the wanted conclusion for $R_1 := f(R)$, $\overline{R}_1 := f(2R)$ and $w_1 \equiv 1$.

Now we assume that the conclusion is true for $n\ge 1$ and consider the case for $n+1$. We set
\begin{equation}\label{eq-different-setting}
\overline{\overline{R}}_i := f(2\overline{R}_i),\;\; \overline{R}'_i :=  f(\overline{R}_i), 1\le i \le n.
\end{equation}
Functional calculus with the conditions (2) and (3) tells us that
\begin{equation}\label{eq-update1}
\overline{R}'_i\overline{R}'_j = \overline{R}'_j\overline{R}'_i = 0,\;\; R_i\overline{R}'_i = \overline{R}'_iR_i = R_i,\;\; 1\le i\ne j\le n.
\end{equation}
Moreover, we can easily see that $f(t)f(2t)=f(t)$, $0\le t\le 1$, so that we have
\begin{equation}\label{eq-update2}
\overline{R}'_i \overline{\overline{R}}_i =  \overline{\overline{R}}_i \overline{R}'_i,\;\; 1\le i \le n.
\end{equation}
Since $\sum^n_{i=1}(\overline{\overline{R}}_i)^2 \in C^*_r(G,\sigma)$ and $G$ is non-discrete we can apply Lemma \ref{lem-cpt} to get a unit vector $h_n \in L^2(G)$ such that
\begin{equation}\label{eq-double-line1}
|\la \left(\Big[\sum^n_{i=1}(\overline{\overline{R}}_i)^2\Big]\cdot u_{n+1}\right) h_n, h_n \ra| < \Big(\frac{\eps^2}{64n}\Big)^2
\end{equation}
and
\begin{equation}\label{eq-double-line2}
|\la \lambda_\sigma(x_{n+1})^*\Big[\sum^n_{i=1}(\overline{\overline{R}}_i)^2\Big] \lambda_\sigma(s_{n+1}) h_n, h_n \ra| < \Big(\frac{\eps^2}{64n}\Big)^2.
\end{equation}
Now we set $w_{n+1}(\cdot) = \la \lambda(\cdot)h_n, h_n \ra \in \mathcal{P}_1(G)$. By Corollary \ref{cor-sigma-pd} we know that $u_{n+1}w_{n+1} \in \mathcal{P}_1(G,\sigma)$ with compact support, so that we can apply Kaplansky density theorem again to get $T\in C^*_r(G,\sigma)$ such that $0\le T \le 1$ and
	\begin{equation}\label{eq-step-n+1-1}
	\la T, u_{n+1}w_{n+1} \ra > 1 - \frac{\eps^2}{8}
	\end{equation}
and
	\begin{equation}\label{eq-step-n+1-2}
	\la\lambda_\sigma(s_{n+1})^* T \lambda_\sigma(s_{n+1}), u_{n+1}w_{n+1} \ra > 1 - \frac{\eps^2}{8}.
	\end{equation}	
We further define
	$$R' := \Big(\prod^n_{i=1}(1-\overline{\overline{R}}_i)\Big)T \Big(\prod^n_{i=1}(1-\overline{\overline{R}}_i)\Big),\;\; R_{n+1} := f(R'),\;\; \overline{R}_{n+1} := f(2R').$$
From \eqref{eq-update2} we can see that $R'\overline{R}'_i = 0$, so that $\overline{R}_{n+1}\overline{R}'_i = 0$ for $1\le i \le n$. Similarly, we also get $\overline{R}'_i\overline{R}_{n+1} = 0$, $1\le i \le n$. Together with \eqref{eq-update1} we now see that the conditions (1), (2) and (3) are satisfied for the sequences $(R_1, \cdots, R_{n+1})$ and $(\overline{R}'_1, \cdots, \overline{R}'_n, \overline{R}_{n+1})$.

Note that the condition \eqref{eq-double-line1} is the same as $|\la \sum^n_{i=1}(\overline{\overline{R}}_i)^2, u_{n+1}w_{n+1} \ra| < \Big(\frac{\eps^2}{64n}\Big)^2$, so that we can apply Lemma \ref{lem-inner-control} (2) together with \eqref{eq-step-n+1-1} to get $\la R', u_{n+1}w_{n+1} \ra>1-\frac{\eps^2}{8}$. By repeating the same argument with \eqref{eq-double-line2}, \eqref{eq-step-n+1-2} and unitary conjugation with respect to $\lambda_\sigma(s_{n+1})$ (i.e. replacing $X$ with $\lambda_\sigma(s_{n+1})^*X\lambda_\sigma(s_{n+1})$) we also get
	$$\la \lambda_\sigma(s_{n+1})^*R'\lambda_\sigma(s_{n+1}), u_{n+1}w_{n+1} \ra>1-\frac{\eps^2}{8}.$$
Finally we appeal to Lemma \ref{lem-inner-control} (1) again to get
	$$\langle f(R') \lambda_\sigma(s_{n+1}) f(R') , \lambda_\sigma(s_{n+1})\bullet u_{n+1}w_{n+1} \rangle >1-\eps,$$
which gives us the condition (4) and the induction procedure is now completed.

For the last statement we set $S_i = R_i\lambda_\sigma(s_i)R_i$ and note that $(S_i)^n_{i=1}$ (respectively, $(S^*_i)^n_{i=1}$) have orthogonal ranges from the conditions (2) and (3). Then, we have
$\displaystyle \|\sum^n_{j=1}S_j\|\le 1$
and
\begin{align}\label{eq-lower-estimate}
\lefteqn{\|\sum^n_{i=1}\mu_i \lambda_\sigma(s_i) \bullet (u_iw_i)\|_{A(G,\sigma)}}\nonumber\\
& \ge |\la \sum^n_{j=1}S_j,  \sum^n_{i=1}\mu_i \lambda_\sigma(s_i) \bullet (u_iw_i)\ra| = |\sum^n_{i,j=1}\mu_i \la S_j, \lambda_\sigma(s_i) \bullet (u_iw_i)\ra|\nonumber\\
& \ge |\sum^n_{i=1}\mu_i \la S_j, \lambda_\sigma(s_i) \bullet (u_iw_i)\ra| - | \sum^n_{i=1}\mu_i \la \sum_{j\ne i}S_j, \lambda_\sigma(s_i) \bullet (u_iw_i)\ra|.
\end{align}
Note that $(u_iw_i)(\cdot) = \la\lambda_\sigma(\cdot)k_i, k_i\ra$ for some unit vector $k_i\in L^2(G)$, so that we have
	\begin{align*}
		|\la \sum_{j\ne i}S_j, \lambda_\sigma(s_i) \bullet (u_iw_i)\ra|
		& = |\la \sum_{j\ne i}\lambda_\sigma(s_i)^* S_j, u_iw_i\ra|\\
		& = |\la \sum_{j\ne i}\lambda_\sigma(s_i)^* S_j k_i, k_i\ra|\\
		& = |\la k_i, \big(\sum_{j\ne i} S^*_j\big) \lambda_\sigma(s_i)k_i\ra|\\
		& = |\la k_i - S^*_i \lambda_\sigma(s_i)k_i, \big(\sum_{j\ne i} S^*_j\big) \lambda_\sigma(s_i)k_i\ra|\\
		& \le \|k_i - S^*_i \lambda_\sigma(s_i)k_i\|_2 < \sqrt{2\eps}
	\end{align*}
since
\begin{align*}
\|k_i - S^*_i \lambda_\sigma(s_i)k_i\|^2_2
& = \la k_i - S^*_i \lambda_\sigma(s_i)k_i, k_i - S^*_j \lambda_\sigma(s_i)k_i \ra\\
& \le 2 - 2 {\rm Re}\la S^*_i \lambda_\sigma(s_i)k_i, k_i\ra\\
& = 2 - 2 |\la S_i, \lambda_\sigma(s_i) \bullet (u_iw_i)\ra| < 2\eps.
\end{align*}
Note that the above estimate is still valid when we replace $(S_i)^n_{i=1}$ with $(z_iS_i)^n_{i=1}$ for any $(z_i)^n_{i=1} \subseteq \tor$, so that \eqref{eq-lower-estimate} gives us the conclusion we wanted.
\end{proof}

\begin{rem}\label{rem-comments-proof-Losert}
Proposition \ref{prop-ell_1} is a corrected and updated version of \cite[Proposition 1]{Losert} to the twisted setting. We made the induction in the proof of \cite[Proposition 1]{Losert} explicit, which was crucial but slightly flawed as originally presented. Precisely speaking, after we replace $\overline{R}_i$ at the end of the proof of \cite[Proposition 1]{Losert} with $\overline{R}'_i := (1-f(4R'))\overline{R}_i(1-f(4R'))$, $1\le i \le n$, where $R_i$ and $R'$ are the operators from \cite[Proposition 1]{Losert}, the desired properties are not guarranteed since we do not know whether $R'$ and $\overline{R}_i $, $1\le i \le n$, are commuting.

In the beginning stage of this project the authors had a different solution for the above mentioned error, but the approach in Proposition \ref{prop-ell_1} actually follows Losert's own solution kindly shared with the authors through private communications \cite{Losert-extra1, Losert-extra2}. Note that the proof of \cite[Proposition 5.3.3]{KanLau} also dealing with \cite[Proposition 1]{Losert} is still incomplete as it is.
\end{rem}

({\it Proof of Theorem \ref{eq-imply-amen}: for non-discrete $G$})

We begin with the equality $M_{cb}(A(G), A(G,\sigma)) = B(G,\sigma)$ and will deduce the hypothesis of Proposition \ref{prop-LUC-mean}, i.e. for any points $s_1, \cdots, s_N \in G$ we would like to construct a mean $M$ on $LUC(G)$ such that $\sum^N_{k=1}d(M,s_k)\le cN$ for some universal constant constant $c<2$ independent of $N$.

{\it Step 1.} For an arbitrary $\eps>0$ and a compact neighborhood $V$ of $e\in G$, we can choose neighborhoods $W$ and $U$ of $e\in G$ such that
	\begin{equation}\label{eq-UW}
	UW\subseteq V
	\end{equation}
and
	\begin{equation}\label{eq-VW-eps}
	\mu(V) < (1+\eps)\mu(W).
	\end{equation}	
By a standard argument we can choose a symmetric neighborhood $U' = (U')^{-1}$ of $e\in G$ such that $U' \subseteq U$ and a function $\varphi \in A(G)$ such that
	\begin{equation}\label{eq-varphi}
	\begin{cases}{\rm supp}\,\varphi\subseteq U,\\ \varphi|_{U'} \equiv 1,\\ \|\varphi\|_{A(G)} \le 1+\eps.\end{cases}
	\end{equation}
Now we choose a function $u \in \mathcal{P}_1(G,\sigma)$ with ${\rm supp}\,u\subseteq U'$. Indeed, we pick a symmetric neighborhood $U''$ of $e\in G$ such that $U'' \subseteq U'$ and $U''U'' \subseteq U'$. Then, the following choice of $u$ is what we wanted:
	\begin{equation}\label{eq-u-choice}
	u(\cdot) := \la \lambda_\sigma(\cdot) \xi, \xi\ra,
	\end{equation}
where $\xi = \sqrt{\mu(U'')}1_{U''} \in L^2(G)$.
	
We finally claim in this step that there is a positive norm 1 element $v = [v_{ij}] \in M_n(VN(G))_*$ and a universal constant $C>0$ such that
	\begin{equation}\label{eq-v-choice}
	\sum^N_{k=1}\|[u \cdot (\lambda(s_k)^*v_{ij})]\|_{M_n(VN(G,\sigma))_*} > CN.
	\end{equation}
For the choice of such $v$ we apply Proposition \ref{prop-ell_1} for $u_k = u$, $1\le k \le N$, to get $(w_k)^N_{k=1}\subseteq \mathcal{P}_1(G)$ such that
	$\|\sum^N_{k=1}\lambda_\sigma(s_k)\bullet (uw_k)\|_{A(G,\sigma)} > N/2.$
This means that
	$$\|\sum^N_{k=1}\lambda_\sigma(s_k)\bullet (uw_k)\|_{M_{cb}(A(G), A(G,\sigma))} > CN$$
for some universal constant $C>0$ coming from the equality $M_{cb}(A(G), A(G,\sigma)) = B(G,\sigma)$, so that we can find a contractive element $v = [v_{ij}] \in M_n(VN(G))_*$ such that
	\begin{align*}
	\lefteqn{\sum^N_{k=1}\|[\big(\lambda_\sigma(s_k)\bullet (uw_k)\big)v_{ij}]\|_{M_n(VN(G,\sigma))_*}}\\
	& \ge \|[\Big(\sum^N_{k=1}\lambda_\sigma(s_k)\bullet (uw_k)\Big)v_{ij}]\|_{M_n(VN(G,\sigma))_*}\\
	& > CN.
	\end{align*}
Since $v$ is a linear combination of four positive contractive elements in $M_n(VN(G))_*$, we may assume that $v$ itself is positive and by multiplying a suitable positive constant we may further assume that $v$ has norm 1.

A direct computation tells us that 
	$$\lambda_\sigma(s^{-1}_k)\bullet \Big[\big(\lambda_\sigma(s_k)\bullet (uw_k)\big)v_{ij}\Big] = \overline{\sigma(s^{-1}_k, s_k)}uw_k \cdot (\lambda(s_k)^*v_{ij}).$$
From (3) of Proposition \ref{prop-two-actions} and the fact that $f \in A(G,\sigma) \mapsto fw_k \in A(G,\sigma)$ is a complete contraction by Proposition \ref{prop-bimodule} and  the condition $w_k \in \mathcal{P}_1(G)\subseteq B(G)$ we get
	\begin{align*}
		\sum^N_{k=1}\|[u \cdot (\lambda(s_k)^*v_{ij})]\|_{M_n(VN(G,\sigma))_*}
		& \ge \sum^N_{k=1}\|[uw_k \cdot (\lambda(s_k)^*v_{ij})]\|_{M_n(VN(G,\sigma))_*} \nonumber\\
		& = \sum^N_{k=1}\|[\big(\lambda_\sigma(s_k)\bullet (uw_k)\big)v_{ij}]\|_{M_n(VN(G,\sigma))_*}\nonumber \\
		& > CN.
	\end{align*}

\vspace{0.3cm}	
{\it Step 2.} We would like to get an upper estimate for $C_k := \|[u \cdot (\lambda(s_k)^*v_{ij})]\|_{M_n(VN(G,\sigma))_*}$, $1 \le k \le N$. We first pick a contractive element $T^k = [T^k_{ij}] \in M_n(VN(G,\sigma))$ such that
	\begin{align*}
		C_k
		& = |\la T^k, [u \cdot (\lambda(s_k)^*v_{ij})] \ra|
		= |\la T^k, [\varphi u \cdot (\lambda(s_k)^*v_{ij})] \ra|\\
		& = |\sum^n_{i,j=1}\la T^k_{ij}, \varphi u \cdot (\lambda(s_k)^*v_{ij}) \ra|
		= |\sum^n_{i,j=1}\la T^k_{ij}\cdot \varphi, u \cdot (\lambda(s_k)^*v_{ij}) \ra|\\
		& = |\sum^n_{i,j=1}\la \Gamma^\sigma(T^k_{ij}\cdot \varphi), u \otimes \lambda(s_k)^*v_{ij} \ra|,
	\end{align*}
where $\Gamma^\sigma$ is the twisted co-multiplication from \eqref{eq-twisted-co-product}. Note that we used the fact that $u = \varphi u$ from \eqref{eq-varphi} for the second equality and (1) of Proposition \ref{prop-two-actions} for the fourth equality. Since $v$ is a positive contractive element of $M_n(VN(G))_*$, we have
	\begin{equation}\label{eq-v-standard}
	v = [v_{ij}] = \omega_{\eta, \eta}
	\end{equation}	
for some unit vector $\eta = [\eta_{ij}] \in H_n \otimes_2 L^2(G) \cong L^2(G;H_n)$, where $H_n$ is the Hilbert space of all $n\times n$ matrices with the canonical inner product $\la X,Y\ra = {\rm Tr}(XY^*)$. Here, we are using the fact that the von Neumann algebras $VN(G)$ and $M_n(VN(G))$ are both in standard form.	Together with \eqref{eq-u-choice} we have
	\begin{align*}
	\sum^n_{i,j=1}\langle \Gamma^\sigma \big(T^k_{ij} \cdot \varphi \big) , u\otimes \lambda(s_k)^*v_{ij}  \rangle
	& = \big \langle\sum_{i,l}\big(\sum_j \Gamma^\sigma(T^k_{ij} \cdot \varphi) (\xi\otimes \eta_{jl})\big), \xi\otimes \lambda(s_k)^*\eta_{il}  \big \rangle\\
	& = \mu(W)^{-1} \cdot I,
	\end{align*}
where $I$ is the integral
	$$I = \int_W \big| \int_{G\times G} \sum_{i,j,l} \Gamma^\sigma(T^k_{ij} \cdot \varphi) (\xi\otimes \eta_{jl})(x,zy)\,\overline{\big(\xi\otimes \lambda(x_k)^*\eta_{il}\big) (x,zy) }\,dx dy  \big| dz.$$
Moreover, we have
	\begin{align*}
	\lefteqn{I \le \int_G \int_{G\times W} \sum_{i,l} \big| \sum_j \Gamma^\sigma(T^k_{ij} \cdot \varphi) (\xi\otimes \eta_{jl})(x,zy) \big| \cdot \big| \big(\xi\otimes \lambda(s_k)^*\eta_{il}\big) (x,zy)  \big| \,dx   dz\,dy} \\
	& = \int_G \int_{G\times Wy} \sum_{i,l} \big| \sum_j \Gamma^\sigma(T^k_{ij} \cdot \varphi) (\xi\otimes \eta_{jl})(x,z) \big| \cdot \big| \big(\xi\otimes \lambda(s_k)^*\eta_{il}\big) (x,z)  \big| \,dx   dz\,\frac{dy}{\Delta_G(y)}\\
	& = \int_G \int_{G\times Wy} \sum_{i,l} \big| \sum_j \Gamma^\sigma(T^k_{ij} \cdot \varphi) (\xi\otimes \eta_{jl}|_{Vy})(x,z) \big| \cdot \big| \big(\xi\otimes \lambda(s_k)^*\eta_{il}\big) (x,z)  \big| \,dx   dz\,\frac{dy}{\Delta_G(y)}\\
	& \le (1+\eps)\int_G \|[\xi \otimes \eta_{ij}|_{Vy}]\|_{L^2(G\times G;H_n)} \cdot  \|[\xi \otimes (\lambda(s_k)^*\eta_{ij})|_{Wy}]\|_{L^2(G\times G;H_n)}   \,\frac{dy}{\Delta_G(y)}\\
	& = (1+\eps)\int_G \|[\eta_{ij}|_{Vy}]\|_{L^2(G;H_n)} \cdot  \|[(\lambda(s_k)^*\eta_{ij})|_{Wy}]\|_{L^2(G;H_n)}   \,\frac{dy}{\Delta_G(y)}.
	\end{align*}
In the above the second equality is from the fact that $({\rm supp}\varphi) Wy\subseteq UWy \subseteq Vy$ and that each $T^k_{ij}\cdot \varphi$ is a $\sigma$-SOT(strong operator topology) limit of the elements of the form $\sum^m_{i=1}\alpha_i\varphi(y_i)\lambda_\sigma(y_i)$ with $(y_i)^m_{i=1}\subseteq {\rm supp}\varphi$, so that $\Gamma^\sigma(T^k_{ij}\cdot \varphi)$ is a $\sigma$-SOT limit of the elements of the form $\sum^m_{i=1}\alpha_i\varphi(y_i)\lambda_\sigma(y_i)\otimes \lambda(y_i)$. The last inequality in the above is by the Cauchy-Schwarz inequality and the fact that $\|[\Gamma^\sigma(T^k_{ij}\cdot \varphi)]\| \le \|[T^k_{ij}]\| \cdot \|\varphi\|_{A(G)} \le 1+\eps$ from \eqref{eq-varphi}.

Finally, we set
	$$a(y) := \frac{\|[\eta_{ij}|_{Vy}]\|_{L^2(G;H_n)}}{\sqrt{\Delta_G(y)}}\;\; \text{and}\;\; b(y) := \frac{\|[(\lambda(x_k)^*\eta_{ij})|_{Wy}]\|_{L^2(G;H_n)}}{\sqrt{\Delta_G(y)}},\; y\in G.$$
to get the estimate
	\begin{equation}\label{eq-C_k-estimate}
	C_k \le \frac{1+\eps}{\mu(W)}\int_Ga(y)b(y)d\mu(y).
	\end{equation}	
Note that it is straightforward to check
	\begin{equation}\label{eq-a-b-norm}
	\int_G a(y)^2d\mu(y) = \mu(V)\;\; \text{and}\;\; \int_G b(y)^2d\mu(y) = \mu(W).
	\end{equation}	

\vspace{0.3cm}	
{\it Step 3.} Finally we define a mean $M = M_V$ on $LUC(G)$ by
	$$M(f) := \int_G f(y) \sum^n_{i,j=1}|\eta_{ij}(y)|^2d\mu(y),\;\; f\in LUC(G).$$
From \eqref{eq-a-b-norm} we have
\begin{align*}
\lefteqn{\Big|   M(f) -        \mu(V)^{-1}   \int_G f(y)a(y)^2 \, dy\Big|}\\
& = \Big|   M(f) -        \mu(V)^{-1}   \int_G f(y) \int_V  \sum_{i,j}  |\eta_{ij}(zy)|^2 \, dz  \, dy\Big|\\
& = \Big|   M(f) -         \mu(V)^{-1}   \int_V \int_G     f(z^{-1}  y)       \sum_{i,j}  |\eta_{ij}(y)|^2 \, dy \, dz \Big|  \\
& \le \sup_{z\in V}  \|  \lambda(z) f-f\|_\infty \,   .
\end{align*}
Similary we get
	$$\Big|   M(\lambda(x_k)f) - \mu(W)^{-1}   \int_G f(y)b(y)^2 \, dy\Big| \le \sup_{z\in W}  \|  \lambda(z) f-f\|_\infty \le \sup_{z\in V}  \|  \lambda(z) f-f\|_\infty,$$
so that we have
	\begin{align*}
	\Big|   M(\lambda(x_k)f - f)\Big|
	& \le 2 \sup_{z\in V}  \|  \lambda(z) f-f\|_\infty\\
	& \;\;\;\; + \Big| \mu(V)^{-1}   \int_G f(y)a(y)^2\,dy - \mu(W)^{-1}   \int_G f(y)b(y)^2 \, dy \Big|\\
	& \le 2 \sup_{z\in V}  \|  \lambda(z) f-f\|_\infty + \mu(W)^{-1}\int_G |f(y)|\cdot | \frac{\mu(W)}{\mu(V)}a(y)^2 - b(y)^2 |dy\\
	& \le 2 \sup_{z\in V}  \|  \lambda(z) f-f\|_\infty + D \cdot \|f\|_\infty,
	\end{align*}
where $D = \eps + \mu(W)^{-1}\int_G |a(y)^2 - b(y)^2 |dy$. By the estimate \eqref{eq-C_k-estimate} and \cite[Lemma 4]{Losert} we have
	\begin{align*}
	D& \le \eps + \mu(W)^{-1}\big((\mu(V)+\mu(W))^2 - \frac{4\mu(W)^2C^2_k}{(1+\eps)^2}\big)^{\frac{1}{2}}\\
	& \le \eps + \big[(2+\eps)^2 - \frac{4C^2_k}{(1+\eps)^2}\big]^{\frac{1}{2}}\\
	& \le 2+3\eps - C^2_k.
	\end{align*}
By summation over $k$	we deduce from \eqref{eq-v-choice} that
	$$\sum^N_{k=1} \Big|   M(\lambda(s_k)f - f)\Big| \le 2N \sup_{z\in V}  \|  \lambda(z) f-f\|_\infty + \|f\|_\infty \big( (2+3\eps)N - NC^2\big).$$
Now, we fix $\eps>0$ such that $\tilde{C} = 2+3\eps - C^2 <2$. Then, we have a family of means $(M_V)$ on $LUC(G)$ satisfying
	$$\sum^N_{k=1} \Big|   M_V(\lambda(s_k)f - f)\Big| \le 2N \sup_{z\in V}  \|  \lambda(z) f-f\|_\infty + \tilde{C}N\|f\|_\infty.$$
By considering a weak$^*$-cluster point in $LUC(G)^*$ we actually get a mean $m$ satisfying
	$$\sum^N_{k=1} \Big|   m(\lambda(s_k)f - f)\Big| \le \tilde{C}N\|f\|_\infty,\;\; f\in LUC(G),$$
which	 forces $G$ to be amenable by Proposition \ref{prop-LUC-mean}.

\begin{rem}\label{rem-comment-cb}
The above proof combined Lemma 3 and Lemma 5 and the proof of Theorem 1 in p. 352-353 of \cite{Losert}. For the extension to the cb-multiplier version we followed the approach of an unpublished note by Ruan \cite{Ruan-unpub}.
\end{rem}

\section*{Acknowledgements} We, the authors, are grateful to Victor Losert for his patient answers to numerous questions of the authors regarding his own solution to the error in \cite[Proposition 1]{Losert}, to Zhong-Jin Ruan for kindly sharing his unpublished manuscript \cite{Ruan-unpub}, to Matthew Daws for kindly communicating to us the standard form of twisted group von Neumann algebras, and to Anthony To-Ming Lau, Nico Spronk for valuable discussions and information.
X. Xiong is supported by the National Natural Science Foundation of China, No. 12371138 and No. W2441002.

\end{document}